\documentclass[]{article}

\usepackage{amsfonts,amsmath,amssymb}
\usepackage{amsthm}
\usepackage[cp866]{inputenc}
\usepackage[english]{babel}
\usepackage[abs]{overpic}

\allowdisplaybreaks \theoremstyle{definition}

\theoremstyle{plain}
\newtheorem{Theorem}{Theorem}
\newtheorem{Corollary}{Corollary}
\newtheorem{Proposition}{Proposition}

\newtheorem{Lemma}{Lemma}

\newcommand{\bfSig}{\mathbf{\Sigma}}
\newcommand{\bfPi}{\mathbf{\Pi}}
\newcommand{\bfDe}{\mathbf{\Delta}}

\newcommand{\bi}{\mathbf{i}}

\title{Wadge Degrees of Classes of $\omega$-Regular $k$-Partitions%\footnote{This paper completes the conference papers \cite{s11,s20} by providing full details for technically involved proofs that were only sketched, and by developing a general approach to other similar problems.}
}
\author{Victor  Selivanov\thanks{This work was supported by the Russian Science Foundation, project 18-11-00100.}
 \\A.P. Ershov Institute of
Informatics Systems SB RAS, Russia\\and\\
Department of Mathematics and Computer Science,\\
	St.~Petersburg State University,\\
	7/9 Universitetskaya nab., Saint Petersburg 199034, Russia
\\{\tt vseliv@iis.nsk.su}
}

\begin{document}
\large
\date{}

 \maketitle

\begin{abstract}
We  develop  a theory of  $k$-partitions of the set of infinite words recognizable by classes of finite automata. The theory enables to complete proofs of existing results about topological classifications of the (aperiodic) $\omega$-regular $k$-partitions, and provides   tools for dealing with other similar questions.  In particular, we characterise the structure of Wadge  degrees of (aperiodic) $\omega$-regular $k$-partitions, prove the decidability of many related problems,  and discuss their complexity.

{\em Keywords.} Wadge reducibility, regular $k$-partition, acceptor, transducer, determinacy, iterated labeled tree, fine hierarchy. 
%{\bf MSC}: Primary   03D15, 03D78, 58J45; Secondary 65M06, 65M25.
\end{abstract}

\section{Introduction}\label{in}

Working in descriptive set theory, W. Wadge \cite{wad84} has shown that the degree structure of Borel sets of $\omega$-words over any finite non-unary alphabet under the many-one reducibility by continuous functions is semi-well-ordered (i.e., it is well founded and has no 3 pairwise incomparable elements).
Working in automata theory independently of W. Wadge, K. Wagner \cite{wag79} has shown that the structure of  regular $\omega$-languages under the continuous reducibility is semi-well-ordered with the corresponding ordinal $\omega^\omega=sup\{\omega, \omega^2,\omega^3,\ldots\}$. Working in computability theory independently of W. Wadge and K. Wagner, the author \cite{s83} discovered a semi-well-ordered structure of ``natural'' $m$-degrees with the corresponding ordinal $\varepsilon_0=sup\{\omega, \omega^\omega,\omega^{\omega^\omega},\ldots\}$. In \cite{s89} (see also \cite{s95,s08a}) we characterised the initial segments of the structure in \cite{s83} by set operations which led to the so called fine hierarchy (FH) of arithmetical sets turned out to be a finitary effective version of the Wadge hierarchy.

In \cite{s98} the Wagner hierarchy  was related to the Wadge hierarchy  and to the FH (see also \cite{cp97,cp99,dr06} for an alternative approach), revealing new connections between descriptive set theory, computability, and automata theory. Later, some results from \cite{wag79,s98} were extended to languages recognized by  more complicated computing devices (see e.g. \cite{d03,s03,fin06} and references therein for an extensive study of, in particular, context-free $\omega$-languages). In this wider context, some important properties of the Wagner hierarchy (e.g., the decidability of levels) usually fail.
It is also natural to investigate variants of the Wagner hierarchy for popular subclasses of regular languages the most important of which is certainly the class of regular aperiodic $\omega$-languages (for brevity, just aperiodic sets). In \cite{s08} a complete aperiodic analogue of the Wagner hierarchy was developed that has its own flavour.

In \cite{s11} we extended the Wagner theory from the regular  sets to the regular $k$-partitions $A:X^\omega\rightarrow \{0,\ldots,k-1\}$ of the set $X^\omega$ of $\omega$-words over a finite alphabet $X$ that essentially coincide with the $k$-tuples $(A_0,\ldots,A_{k-1})$ of pairwise disjoint regular sets
satisfying $A_0\cup\cdots\cup A_{k-1}=X^\omega$ (note that the
$\omega$-languages are in a bijective correspondence with the
2-partitions of $X^\omega$).  Motivations for this generalization
come from the fact that similar objects were studied e.g. in computability theory \cite{s83}, descriptive set theory \cite{he93},  and complexity
theory \cite{ko00}. In \cite{s20} we sketched a Wagner theory for the $\omega$-regular aperiodic $k$-partitions in parallel to the theory in \cite{s11}. 
Note that the extension from sets to $k$-partitions for $k>2$ is  non-trivial. It required to develop a machinery of iterated labeled trees and of the FH of $k$-partitions (systematized in \cite{s12}) turned out crucial for the subsequent partial extension of the Wadge theory to $k$-partitions  \cite{s17} and, as a  concluding step, to the   Borel $Q$-partitions for arbitrary better quasiorder $Q$ \cite{km19}. 

An important advantage of the Wagner hierarchy over  the Wadge hierarchy and the aforementioned hierarchies of $\omega$-context-free languages is the decidability of all natural associated algorithmic problems. Moreover, many such problems for hierarchies of sets are PTIME-computable \cite{kpb,wy}. In a sense, the Wagner theory may be considered as a constructive version of a small fragment of the Wadge theory. Thus, the Wagner hierarchy is important not only as a fundamental classification of $\omega$-regular languages but also as a starting point in identifying a constructive part of the rather non-constructive Wadge hierarchy.  

This paper completes the conference papers \cite{s11,s20} by providing full details for technically involved proofs that were only sketched. It also unifies the proofs and develops a general technique that might be of use in treating similar questions for other classes of languages and $k$-partitions. We illustrate this by discussing classes of $\omega$-regular $k$-partitions related to the languages of finite  words studied in \cite{str,s09,rt,ch12}. In fact, the Wadge degrees of $k$-partitions in some of these classes may be the same as for the regular $k$-partitions but their ``automatic'' versions are usually more subtle.

To make the paper self-contained, we include in the next two sections some older material. In Section \ref{prel} we recall relevant  facts about regular acceptors and transducers, and provide some of their analogues for the aperiodic regular languages and $k$-partitions, and of the classical B\"uchi-Landweber theorem about regular Gale-Stewart games. In Section \ref{fine} we recall related facts about the iterated labelled trees and FHs, and give some new information; altogether, these facts provide  useful tools for dealing with hierarchies of  $k$-partitions of $\omega$-words. In Section \ref{regular} we prove our main results about the FH of (aperiodic) $\omega$-regular $k$-partitions in a uniform way, including ``automatic'' versions of the Wadge reducibility. In Section \ref{compl} we discuss computability and complexity of some related algorithmic problems. We show that the decidability of Wagner hierarchy survives under the extensions from sets to $k$-partitions for $k>2$ and discuss complications that arise when one attempts to extend also the PTIME-decidability from \cite{kpb,wy}. In Section \ref{other} we discuss possible variants of our results for some other classes of  $k$-partitions. 

We use standard set-theoretic notation. For  sets $A$ and $S$, $P(S)$ is the class of subsets of $S$ and $S^A$ is the class of all functions from $A$ to $S$. For a class ${\cal C}\subseteq P(S)$, $\check{\cal C}$ is the dual class $\{S\setminus C\mid C\in{\mathcal
C}\}$ ($S\setminus C$ is often denoted by $\overline{C}$), and $BC({\cal C})$ is the Boolean closure of ${\cal C}$. For a subalgebra $\mathcal{B}$ of $(P(S),\cup,\cap,\bar{\;},\emptyset,S)$ and an integer $k>1$, let $\mathcal{B}_k=\{B\in k^S\mid\forall i<k(B^{-1}(l)\in\mathcal{B})$ be the set of $k$-partitions of $S$ whose components are in $\mathcal{B}$. We assume familiarity with notions from  logic, including the notions of structure and  quotient-structure.

\section{Acceptors, transducers, reducibilities}\label{prel}

In this section we recall (with some modification and adaptation)  notation, notions and facts about automata used in subsequent sections, and prove some new facts. For additional information see e.g. \cite{pp,th90,th96}.

\subsection{Acceptors and $k$-acceptors}\label{accept}

Fix a finite  alphabet $X$ containing more than one symbol (for
simplicity we may assume that $X=\bar{m}=\{0,\ldots,m-1\}$ for an integer $m>1$, so  $0,1$ are always in $X$). Note that usually we work with the fixed alphabet $X$ but sometimes
we are forced to consider several alphabets simultaneously. The
``fixed-alphabet mode'' is the default one.

Let $X^\ast$, $X^+$, and $X^\omega$ denote
resp. the sets  of all words, all nonempty words, and  all $\omega$-words  over $X$. 
Let $\varepsilon$ be the empty word and $X^{\leq\omega}=X^\ast\cup X^\omega$.
We use   standard notation concerning words and
$\omega$-words. For
$w\in X^\ast$ and $\xi\in X^{\leq\omega}$, $w\sqsubseteq\xi$ means
that $w$ is a substring of $\xi$, $w\cdot\xi=w\xi$ denote the
concatenation, $l=|w|$ is the length of $w=w(0)\cdots w(l-1)$. For
$w\in X^\ast,W\subseteq X^\ast$ and $A\subseteq X^{\leq\omega}$, let
$w\cdot A=\{w\xi :\xi\in A\}$ and $W\cdot A=\{w\xi :w\in W,\xi\in
A\}$. For $k,l<\omega$ and $\xi\in X^{\leq\omega}$, let
$\xi[k,l)=\xi(k)\cdots\xi(l-1)$ and $\xi\upharpoonright_k=\xi[0,k)$.

By an {\em automaton} (over $X$) we mean a triple ${\mathcal M}=(Q,f,in)$
consisting  of a finite non-empty set $Q$ of states,  a transition
function $f:Q\times X\rightarrow Q$ and an initial state $in\in Q$.
The  function $f$ is  extended to the function
$f:Q\times X^*\rightarrow Q$  by induction
$f(q,\varepsilon)=q$ and $f(q,u\cdot x)=f(f(q,u),x)$, where $u\in
X^*$ and $x\in X$. Similarly, we may define the function $f:Q\times
X^\omega\rightarrow Q^\omega$ by $f(q,\xi)(n)=f(q,\xi\upharpoonright_n)$. Instead of $f(q,u)$ we often write $q\cdot u$.

Associate
with any automaton ${\mathcal M}$ the set of {\em cycles} $C_{\mathcal M}=\{f_{\mathcal
M}(\xi)\mid\xi\in X^\omega\}$ where $f_{\mathcal M}(\xi)$ is the set of
states that occur infinitely often in the sequence $f(in,\xi)\in
Q^\omega$.
A {\em Muller acceptor} is a pair $({\mathcal M},{\mathcal F})$ where ${\mathcal M}$ is an
automaton and ${\mathcal F}\subseteq C_{\mathcal M}$; it recognizes the set
$ L({\mathcal M},{\mathcal F})=\{\xi\in X^\omega\mid f_{\mathcal M}(\xi)\in{\mathcal
F}\}$. The Muller  acceptors recognize exactly the
{\em regular $\omega$-languages}. Let $\mathcal{R}$ denote the set of all such languages; this set is closed under the Boolean operations. 
 
An automaton ${\cal M}=(Q,X,f)$ is {\em aperiodic} if for all $q\in Q$, $u\in
X^+$ and $n>0$ the equality $f(q,u^n)=q$ implies $f(q,u)=q$. This is
equivalent to saying that for all $q\in Q$ and $u\in X^+$ there is
$m<\omega$ with $f(q,u^{m+1})=f(q,u^m)$. An acceptor is aperiodic if so is the corresponding
automaton. A  language $L\subseteq X^*$ ($L\subseteq X^\omega$) is {\em aperiodic} if it is recognized by an aperiodic (Muller) acceptor. Let $\mathcal{A}$ denote the set of all such $\omega$-languages; this set is closed under the Boolean operations. The aperiodic sets are precisely those which satisfy a fixed first-order sentence (see also Section \ref{other}). 

A regular $k$-partition $L$ may be specified by a $k$-tuple of Muller acceptors that recognize the components $L_0,\ldots,L_{k-1}$ but for our purposes we need a slightly different presentation introduced in \cite{s07a}. 
An (aperiodic) {\em Muller  $k$-acceptor} is a pair $({\mathcal M},A)$ where
${\mathcal M}$ is an (aperiodic) automaton and $A: C_{\mathcal M}\rightarrow k$ is a
$k$-partition of $C_{\mathcal M}$. The Muller  $k$-acceptor recognises the aperiodic
$k$-partition $L({\mathcal M},A)=A\circ f_{\mathcal M}$ where $f_{\mathcal
M}:X^\omega\rightarrow C_{\mathcal M}$ is  defined above. 
 
\begin{Proposition}\label{ch}
A $k$-partition $L:X^\omega\rightarrow \bar{k}$ is regular (aperiodic) iff it is
recognised by an (aperiodic) Muller $k$-acceptor.
\end{Proposition}

{\em Proof.} We consider only the non-trivial direction.
Let $L$ be a regular (aperiodic) $k$-partition and $k>2$ (for $k=2$ the assertion is obvious). Then
$L_l$ is regular (aperiodic) for every $l<k$, hence $L_l=L({\mathcal
M}_l,{\mathcal F}_l)$ for some (aperiodic) Muller acceptors $({\mathcal M}_l,{\mathcal
F}_l)$. Let ${\mathcal M}=(Q,f,in)$ be the product of the automata $\mathcal{M}_0,\ldots,\mathcal{M}_{k-2}$ where $Q=Q_0\times\cdots\times
Q_{k-2}$,
$f((q_0,\ldots,q_{k-2}),x)=(f_0(q_0,x),\ldots,f_{k-2}(q_{k-2},x))$
and $in=(in_0,\ldots,in_{k-2})$. By Proposition 4 in \cite{s08}, ${\mathcal M}$ is aperiodic if all $\mathcal{M}_0,\ldots,\mathcal{M}_{k-2}$ are aperiodic. We have
$pr_l(f_{\mathcal M}(\xi))=f_{{\mathcal M}_l}(\xi)$ for all $l<k-1$ and
$\xi\in X^\omega$, where $pr_l:Q\rightarrow Q_l$ is the projection
to the $l$-th coordinate. Since  $L_l$ are
pairwise disjoint, so are also  $pr^{-1}_l({\mathcal F}_l)$. Let
$A: C_{\mathcal M}\rightarrow \bar{k}$ be the unique $k$-partition of $C_{\mathcal M}$
satisfying $A^{-1}(l)=pr^{-1}_l({\cal F}_l)$ for all $l<k-1$. Then the (aperiodic) Muller $k$-acceptor
$({\mathcal M},A)$ recognises $L$.
 \qed

\subsection{Transducers and games}\label{trans}

The set $X^\omega$ carries the {\em Cantor topology} with the open sets
$W\cdot X^\omega$, where $W\subseteq X^\ast$. The {\em Borel sets} in $X^\omega$ are obtained by closing the class of open sets by the operations of complement and countable unions. Let $\mathbf{\Sigma}^0_n,\mathbf{\Pi}^0_n,\mathbf{\Delta}^0_n$ denote levels of the Borel hierarchy in $X^\omega$ \cite{ke}, so, in particular, $\bfSig^0_1$ is the class of open sets, $\bfPi^0_n=\check{\bfSig}^0_n$, ${\mathbf\Delta}^0_n=\bfSig^0_n\cap\bfPi^0_n$, and $\bfSig^0_{n+1}$ is the class of countable unions of $\bfPi^0_n$-sets. It is well known and easy to see that ${\mathcal R}\subset BC({\mathbf\Sigma}^0_2)\subset {\mathbf\Delta}^0_3$.   

The continuous functions on $X^\omega$ are also called here, following \cite{wag79}, {\em continious asynchronous functions}, or $CA$-functions. A {\em
continuous synchronous function}, or just $CS$-function, is a function $f:X^\omega\rightarrow X^\omega$
satisfying $f(\xi)(n)=\phi(\xi\upharpoonright_{(n+1)})$ for some
$\phi:X^\ast\rightarrow X$; in descriptive set theory such functions are known as Lipschitz functions. Clearly, every  $CS$-function is a $CA$-function. Both classes of functions are closed under composition.

A {\em synchronous transducer} (over  $X,Y$) is a tuple
${\mathcal T}=(Q,X,Y,f,g,in)$, also written as ${\mathcal T}=({\mathcal
M},Y,g,in)$, consisting of an automaton ${\mathcal M}$ as above, an
initial state $in$ and an output function $g:Q\times X\rightarrow Y$.
The output function is extended to the function  $g:Q\times
X^*\rightarrow Y^*$ defined by induction
$  g(q,\varepsilon)=\varepsilon$ and $g(q,u\cdot x)=g(q,u)\cdot g(f(q,u),x),$
and to the function  $g:Q\times X^\omega\rightarrow Y^\omega$ defined by
\begin{equation}
g(q,\xi)=g(q,\xi(0))\cdot g(f(q,\xi(0)),\xi(1))\cdot
g(f(q,\xi[0,2)),\xi(2))\cdots.
\end{equation}
In other notation, $g(q,\xi)=lim_ng(q,\xi\upharpoonright_n)$. The transducer
${\mathcal T}$ {\em computes} the function $g_{\mathcal
T}:X^\omega\rightarrow Y^\omega$ defined by $g_{\mathcal
T}(\xi)=g(in,\xi)$. 

 {\em Asynchronous transducers} are defined in the same way, only now the
output function $g$ maps $Q\times X$ into $Y^*$. As a result, the
value $g(q,\xi)$ defined as in (1) is now in $Y^{\leq\omega}$, and
 $g_{\cal T}:X^\omega\to Y^{\leq\omega}$.
Functions computed by  synchronous (resp. asynchronous) transducers are called $DS$-functions (resp. $DA$-funct\-ions). Both  classes of functions are closed under composition \cite{wag79}.
A transducer ${\mathcal T}=({\mathcal
M},Y,g,in)$ is {\em aperiodic} if  ${\mathcal
M}$ is aperiodic. Functions computed by aperiodic synchronous (resp. asynchronous) transducers are called $AS$-functions (resp. $AA$-funct\-ions). By Proposition 10 in \cite{s08}, both  classes of functions are closed under composition. Obviously, every
$AS$-function (resp. $AA$-function) is a $DS$-function (resp. $DA$-function), and every
$DS$-function (resp. $DA$-function) is a $CS$-function (resp. $CA$-function).

We associate with any  $A\subseteq (X\times
Y)^\omega$ the {\em Gale-Stewart game} $G(A)$  played by two players $0$
and $1$ as follows. Player $0$ chooses a letter $x_0\in X$, then
player $1$ chooses a letter $y_0\in Y$, then  $0$ chooses  $x_1\in
X$, then $1$ chooses  $y_1\in Y$ and so on. Each player knows all
the previous moves. After $\omega$ moves, player $0$ (resp. player $1$) has constructed
a word $\xi=x_0x_1\cdots\in X^\omega$ (resp. $\eta=y_0y_1\cdots\in Y^\omega$). Player $1$ wins
if $\xi\times\eta=(x_0,y_0)(x_1,y_1)\cdots\in A$, otherwise player $0$
wins.

A {\em strategy for player} $1$ (player $0$) in the game $G(A)$ is a
function $h:X^+\rightarrow Y$ (respectively, $h:Y^*\rightarrow X$)
that prompts the player $1$'s move  (respectively, the player $0$'s
move) for any finite string of the opponent's previous moves. The strategies for player $1$ (for $0$) are in a
bijective correspondence with the $CS$-functions
$h:X^\omega\rightarrow Y^\omega$ (respectively, with the delayed
$CS$-functions $h:Y^\omega\rightarrow X^\omega$) \cite{s08}; we identify
strategies with the corresponding $CS$-functions.
A strategy $h$ for player $1$ (player $0$)  in the game $G(A)$ is
{\em winning} if the player always wins  following the
strategy, i.e. if  $\xi\times h(\xi)\in A$ for all $\xi\in X^\omega$
(resp. $h(\eta)\times\eta\in \overline{A}$ for all $\eta\in
Y^\omega$). 

One of the best results of descriptive set theory is the Martin
determinacy theorem (see e.g. \cite{ke})  stating that every Borel set
is determined. Note that, since any regular set is Borel, this
implies the determinacy of regular sets. One of the best results of
automata theory is the B\"{u}chi-Landweber regular determinacy
theorem stating that for any regular set $A$ the winner in $G(A)$
may be computed effectively, (s)he has a winning strategy which is a
$DS$-function, and the strategy is also computed effectively. As shown in  Theorem 1 of \cite{s08}, for any aperiodic set $A\subseteq (X\times Y)^\omega$, one of the players has a winning strategy in $G(A)$, the winner is computable and  has an $AS$-winning strategy which is also  computed effectively. Below we refer to the latter result as the {\em aperiodic determinacy theorem}.

\subsection{Reducibilities on $k$-partitions}\label{red0}

Let $\mathcal{F}$ be a set of unary functions on $X^\omega$ that is closed under composition and contains the identity function.
For $A,B\in k^{X^\omega}$, $A$ is  {\em $\mathcal{F}$-reducible} to  $B$ (in symbols $A\leq_\mathcal{F}B$), if $A=B\circ f$ for
some  $f\in\mathcal{F}$. The relation $\leq_\mathcal{F}$ is a preorder on $k^{X^\omega}$, the induced equivalence relation is denoted by $\equiv_\mathcal{F}$; the $\equiv_\mathcal{F}$-equivalence classes are called {\em $\mathcal{F}$-degrees}. For  $\mathcal{C}\subseteq k^{X^\omega}$, $\mathcal{C}$ is {\em closed under $\leq_\mathcal{F}$} if $D\leq_\mathcal{F}C\in\mathcal{C}$ implies $D\in\mathcal{C}$; a $k$-partition $D$ is {\em $\mathcal{C}$-hard}  (in symbols, $\mathcal{C}\leq_\mathcal{F}D$) if $C\leq_\mathcal{F}D$ for  every $C\in\mathcal{C}$; a $k$-partition $C$ is {\em $\mathcal{C}$-complete in $\mathcal{C}$} (in symbols, $C\equiv_\mathcal{F}\mathcal{C}$) if $C\in\mathcal{C}\leq_\mathcal{F}C$.

For the classes of functions introduced in Subsection \ref{trans}, we obtain   reducibilities $\leq_{CA},\leq_{CS},\leq_{DA},\leq_{DS},\leq_{AA},\leq_{AS}$.  Clearly,
 $\leq_{AS}\subseteq\leq_{DS}\subseteq\leq_{CS}\subseteq\leq_{CA}$
  and $\leq_{AA}\subseteq\leq_{DA}\subseteq\leq_{CA}$. From the results in \cite{wag79,s08} it follows that $\mathcal{R}_k$  is closed  under $\leq_{DA},\leq_{DS}$, $\mathcal{A}_k$  is closed  under $\leq_{AA},\leq_{AS}$, but both classes are not closed   under $\leq_{CS}$, $\leq_{CA}$.

%By $\mathbf{\Sigma}^0_2$-function we mean a unary function $f$ on $X^\omega$ such that $f^{-1}(A)\in\mathbf{\Sigma}^0_2$ for every $A\in\mathbf{\Sigma}^0_2$; we denote the corresponding reducibility by $\leq_2$.

We conclude this subsection by a result demonstrating that the non-effective $CS$-reducibility (that is clearly not well suited for automata theory), being restricted to classes recognized by automata, is equivalent to suitable ``automatic'' reducibilities. A similar result also holds for $CA$-reducibility (see Theorem \ref{redcoin} below) but the proof depends on some additional facts.

\begin{Proposition}\label{redcoin1}
 The relation $\leq_{CS}$  coincides with $\leq_{DS}$ on $\mathcal{R}_k$, and with $\leq_{AS}$ on $\mathcal{A}_k$.
\end{Proposition}

{\em Proof.} Both assertions are proved similarly, so we prove only the second one. It suffices to show that, for all $A,B\in\mathcal{A}_k$, $A\leq_{CS}B$ implies $A\leq_{AS}B$.
Let $A\leq_{CS}B$ via a $CS$-function $f:X^\omega\to X^\omega$. Consider the game $G(A,B)$  where players  produce resp. $\xi$ and $\eta$ from $X^\omega$ as   in Subsection \ref{trans}; let player 1 win iff $A(\xi)=B(\eta)$, i.e. $\xi\in A_i\leftrightarrow \eta\in B_i$ for every $i<k$. Then $f$ is a winning strategy for player 1. Since all the components $A_i,B_i$ are aperiodic and $\mathcal{A}$ is closed under the Boolean operations, $G(A,B)$ is  aperiodic. By the aperiodic determinacy, player 1 has an $AS$-winning strategy $g$. Thus, $A\leq_{AS}B$ via $g$.
 \qed

\section{The fine hierarchy of $k$-partitions}\label{fine}

In this section, we briefly recall some  notions and facts about FHs and prove some new facts. Altogether, this gives  a technical tool for proving the main results this paper. For additional details see \cite{s12,s21}.

\subsection{Preorders and semilattices}\label{semilat}

We assume the reader to be familiar with standard terminology and notation related to parially ordered sets (posets) and preorders.  Recall that a {\em semilattice} is a structure $(S;\sqcup)$ with binary operation $\sqcup$ such that $(x\sqcup y)\sqcup z=x\sqcup (y\sqcup z)$, $x\sqcup y= y\sqcup x$  and $x\sqcup x=x$, for all $x,y,z\in S$. By $\leq$ we denote the induced partial order on $S$: $x\leq y$ iff $x\sqcup y=y$. The operation $\sqcup$ can be recovered from $\leq$ since $x\sqcup y$ is the supremum of $x,y$ w.r.t. $\leq$. The semilattice is {\em distributive} if $x\leq y\sqcup z$ implies that $x= y'\sqcup z'$ for some $y'\leq y$ and $z'\leq z$. All semilattices considered in this paper are distributive (sometimes after adjoining a new smallest element denoted by $\bot$). A semilattice $(S;\sqcup,\leq)$ is a {\em d-semilattice} if it becomes distributive after adjoining to $S$ a new smallest element $\bot$.

A non-smallest element $x$ of the semilattice $S$ is {\em join-reducible} if it  can be represented as the  supremum of some elements strictly below $x$. Element $x$  is {\em join-irreducible} if it is not join-reducible.  We denote by $I(S;\sqcup,\leq)$ the set of join-irreducible elements of a semilattice $(S;\sqcup,\leq)$. If $S$ is distributive then $x$  is  join-irreducible iff $x\leq y\sqcup z$ implies that $x\leq y$ or $x\leq z$. By a {\em decomposition of $x$} we mean a representation $x=x_0\sqcup\cdots\sqcup x_n$ where the {\em components $x_i$} are join-irreducible and pairwise incomparable. Such a decomposition is {\em canonical} if it is unique up to a permutation of the components. Clearly, if $S$ 
 is a well founded semilattice then any non-smallest element $x\in S$ has a decomposition, and if $S$ is distributive then $x$ has a canonical decomposition.
 
To simplify notation, we often apply the terminology about posets to preorders meaning the corresponding quotient-poset. Similarly, the term ``semilattice''  will also be applied to structures $(S;\sqcup,\leq)$ where $\leq$ is a preorder on $S$ such the quotient-structure under the induced equivalence relation $\equiv$ is a ``real'' semilattice with the partial order induced by $\leq$ (thus, we avoid  precise but more complicated terms like ``pre-semilattice''). We call preorders (or pre-semilattices) $P,Q$ {\em equivalent} (in symbols, $P\simeq Q$) if their quotient-posets (resp., quotient-semilattices) are isomorphic. For subsets $A,B\subseteq S$ of a preorder $(S;\leq)$ we write $A\equiv B$ if every element of $A$ is equivalent to some element of $B$ and vice versa.
%This terminology is especially convenient in the situation (which is typical below) when one operation $\sqcup$ on $S$ induces the pre-semilattice structure for several preorders on $S$.

We associate with any poset $Q$ the preorder $(Q^*;\leq^*)$ where $Q^*$ is the set of non-empty finite subsets of $Q$, and $S\leq^*R$ iff $\forall s\in S\exists r\in R(s\leq r)$. Let $Q^\sqcup$ be the quotient-poset of $(Q^*;\leq^*)$ and $\sqcup$ be the operation  of supremum in $S$ induced by the operation of  union in $Q^*$. Then $Q^\sqcup$ is a d-semilattice the join-irreducible elements of which coincide with the elements induced  by the singleton sets in $Q^*$ (the new smallest element $\bot$ corresponds to the empty subset of $Q$); thus, $(I(Q^\sqcup);\sqcup,\leq^*)\simeq Q$. Any element of $Q^\sqcup_\bot$ has a canonical decomposition. If $Q$ is well founded then so is also $Q^\sqcup$. The construction $Q\mapsto Q^\sqcup$ is a functor from the category of preorders to the category of semilattices. We will use the following easy fact.

\begin{Proposition}\label{slattice}
Let $f:Q\to I(S)$ be a monotone function from a poset $Q$ to  the set of join-irreducible elements of a semilattice $S$. Then there is a unique semilattice homomorphism $f^\sqcup:Q^\sqcup\to S$ extending $f$. If $f$ is an embedding and $S$ is distributive then $f^\sqcup$ is an embedding.
\end{Proposition}

\subsection{Iterated labeled posets}\label{lab}

Here we discuss iterated labeled posets and forests (introduced in \cite{s11} and systematized in some further publications including \cite{s12}) that are used as  notation systems for the FHs of $k$-partitions. 

Let $(P;\leq)$ be a finite poset; if $\leq$ is clear from the context, we simplify the notation of the poset to $P$. Any subset of $P$ may be considered as a poset with the induced
partial ordering. The {\em rank of a finite poset $P$} is the cardinality of a longest chain in $P$. By a {\em  forest} we mean a finite poset  in which every lower cone $\downarrow{x}$, $x\in P$, is a chain. A {\em  tree} is a forest with the least element (called the {\em  root} of the
tree). 

The ``abstract'' trees (forests) just defined are for almost all purposes equivalent to their isomorphic copies realised as initial segments of $(\omega^*;\sqsubseteq)$ (resp. $(\omega^+;\sqsubseteq)$) where $\sqsubseteq$ is the prefix relation on finite strings of naturals. Below we often work with such ``concrete'' copies which enable to use convenient standard notation for strings.

Let $(Q;\leq)$ be a preorder. A {\em $Q$-poset}  is a triple
$(P,\leq,c)$ consisting of a finite nonempty poset $(P;\leq)$,
$P\subseteq\omega$,  and a labeling $c:P\rightarrow Q$. A {\em morphism}
$f:(P,\leq,c)\rightarrow(P^\prime,\leq^\prime,c^\prime)$ between
$Q$-posets is a monotone function
$f:(P;\leq)\rightarrow(P^\prime;\leq^\prime)$ satisfying $\forall
x\in P(c(x)\leq c^\prime(f(x)))$. 
The {\em $h$-preorder} $\leq_h$ on ${\mathcal P}_Q$ is defined as
follows: $P\leq_h P^\prime$, if there is a  morphism
$f:P\rightarrow P^\prime$. Let ${\mathcal P}_Q$, ${\mathcal
F}_Q$, and ${\mathcal T}_Q$  denote the sets of
all finite $Q$-posets, $Q$-forests, and $Q$-trees, respectively.
For the particular case $Q=\bar{k}=\{0,\cdots,k-1\}$ of
antichain with $k$ elements we denote the corresponding preorders by ${\mathcal
P}_k$, ${\mathcal F}_k$, and ${\mathcal T}_k$. For any $q\in Q$ let $s(q)\in{\mathcal T}_Q$ be the singleton tree labeled by $q$; then $q\leq r$ iff $s(q)\leq_hs(r)$. Identifying $q$ with $s(q)$, we may think that $Q$ is a substructure of ${\mathcal T}_Q$.

The structure $({\mathcal F}_Q;\leq_h,\sqcup)$ is a semilattice equivalent to  $({\mathcal T}^\sqcup_Q;\leq_h,\sqcup)$ above. The supremum operation is given by the disjoint union $F\sqcup G$ of $Q$-forests $F,G$, the join-irreducible elements are precisely the elements $h$-equivalent to trees. In this paper, the  iterations $Q\mapsto{\mathcal T}_{{\mathcal T}_Q}$, $Q\mapsto{\mathcal F}_{{\mathcal T}_Q}$, and $Q\mapsto{\mathcal P}_{{\mathcal P}_Q}$ of these constructions are especially relevant. Using the identification $q=s(q)$, we may think that ${\mathcal T}_Q$ is a substructure of ${\mathcal T}_{{\mathcal T}_Q}$. Define the binary operation  $\cdot$ on ${\mathcal F}_{{\mathcal T}_Q}$ as
follows: $F\cdot G$ is obtained by adjoining a copy of $G$ below any
leaf of $F$. One easily checks that this operation is associative (i.e. $(F\cdot G)\cdot H\equiv_hF\cdot (G\cdot H)$) but not commutative (this was the reason for changing the notation $+$ for this operation in \cite{s11} to  $\cdot $). For $F\in\mathcal{F}_Q$, let $r(F)=\bigsqcup\{c(x)\mid x\in F\}$; then $r:\mathcal{F}_Q\to Q^\sqcup$ is a semilattice homomorphism such that $q=r(s(q))$ for every $q\in Q$.

Recall that a {\em well quasiorder}  (wqo) is a preorder that has
neither infinite descending chains nor infinite antichains.
A famous Kruskal's theorem  implies that if $Q$ is a wqo then $({\mathcal
F}_Q;\leq_h)$ and $({\mathcal T}_Q;\leq_h)$  are wqo's; it is not hard to see that $({\mathcal P}_Q;\leq_h)$ is, in general, not a wqo. Note that the iterated preorders
${\mathcal T}_{{\mathcal T}_Q}$ and ${\mathcal F}_{{\mathcal T}_Q}$
are wqo's whenever $Q$ is a wqo.

Define the sequence $\{\mathcal{T}_k(n)\}_{n<\omega}$ of preorders by induction on $n$
as follows: $\mathcal{T}_k(0)=\overline{k}$  and
$\mathcal{T}_k(n+1)=\mathcal{T}_{\mathcal{T}_k(n)}$. The sets $\mathcal{T}_k(n)$, $n<\omega$, are pairwise disjoint but, identifying
the elements $i$ of $\overline{k}$ with the corresponding singleton trees $s(i)$   labeled by $i$ (which are precisely the
minimal elements  of $\mathcal{T}_k(1)$), we may think that
$\mathcal{T}_k(0)\sqsubseteq\mathcal{T}_k(1)$, i.e. the quotient-poset of the first preorder is an initial segment of the quotient-poset of the other. This also induces an embedding of $\mathcal{T}_k(n)$ into $\mathcal{T}_k(n+1)$ as an initial segment, so (abusing notation) we may think that $\mathcal{T}_k(0)\sqsubseteq\mathcal{T}_k(1)\sqsubseteq\cdots$, hence
$\mathcal{T}_k(\omega)=\bigcup_{n<\omega}\mathcal{T}_k(n)$ is a wqo w.r.t. the induced preorder which we also denote $\leq_h$. We often simplify $\mathcal{T}_k(n)^\sqcup$ to $\mathcal{F}_k(n)$. The embedding $s$ is extended to $\mathcal{T}_k(\omega)$ by defining $s(T)$ as the singleton tree labeled by $T$. Note that ${\mathcal T}_{{\mathcal T}_k}={\mathcal T}_k(2)$ and ${\mathcal F}_{{\mathcal T}_k}={\mathcal F}_k(2)$. Initial segments of $(\mathcal {F}_2(1);\leq_h)$ for $k=2,3$ are depicted below.\footnote{I thank Anton Zhukov for the help with making the pictures.}

\begin{center}
\includegraphics[scale=0.6]{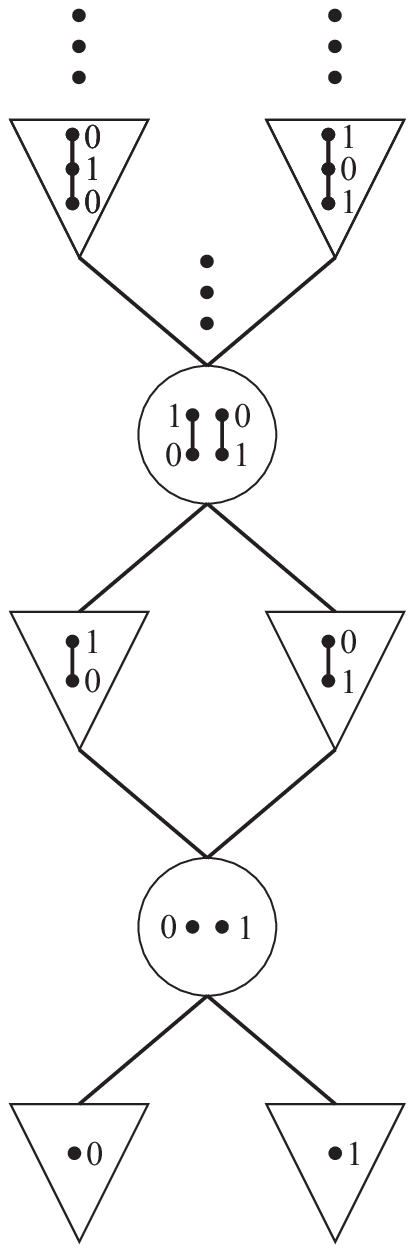}

Fig. 1. An initial segment of $(\mathcal {F}_2(1);\leq_h)$.
\end{center}

\begin{center}
\includegraphics[scale=0.5]{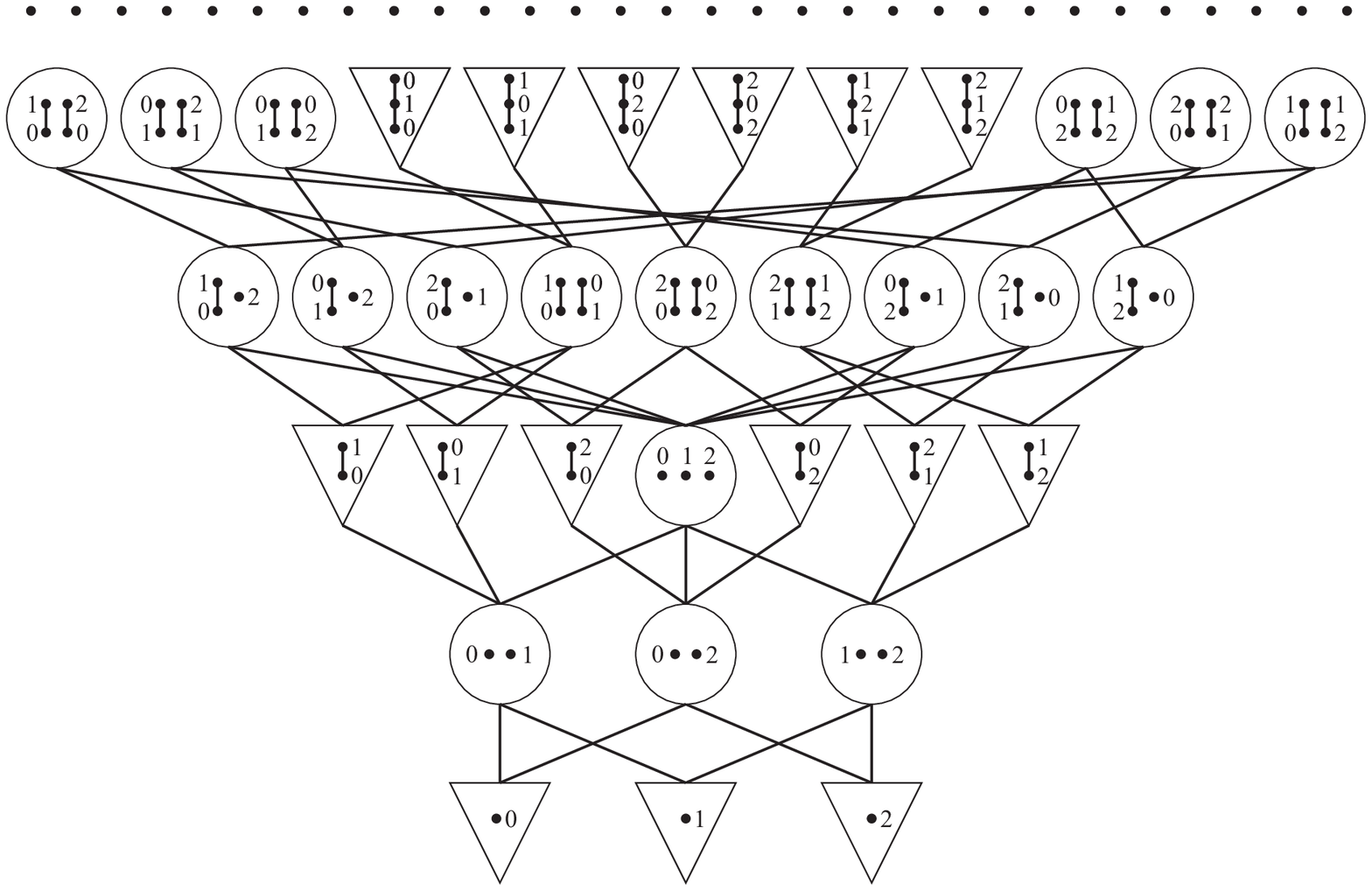}

Fig. 2. An initial segment of $(\mathcal{F}_3(1);\leq_h)$.
\end{center}

Note that while $\mathcal {F}_2(1)$ is semi-well-ordered with rank $\omega$, $\mathcal {F}_k(1)$ for $k>2$ is a wqo of rank $\omega$ having antichains of arbitrary finite size. The whole structure $\mathcal {F}_2(\omega)$  is also semi-well-ordered but with larger rank $\varepsilon_0$ (see  Proposition 8.28 in \cite{s12}). This structure is isomorphic to the FH of arithmetical sets in \cite{s83,s89} mentioned in the Introduction. The triangle levels (induced by trees) correspond to the ``non-self-dual'' $\Sigma$- and $\Pi$-levels of this hierarchy. More precisely, the $\Sigma$-levels (resp. $\Pi$-levels) correspond to (hereditary) 0-rooted (resp., 1-rooted) trees. According to Fig. 2,  the preorder $\mathcal{F}_k(1)$ for $k>2$ is much more complicated than for $k=2$. Nevertheless, the (generalised) non-self-dual levels of the corresponding FHs of $k$-partitions will again correspond to trees (depicted as triangles).

By Proposition 8.7(2) of \cite{s12}, for any finite $Q$-poset $(P,\leq,c)$ there
exist a finite $Q$-forest $F=F(P)$ of the same rank as $P$ (obtained by a bottom-up unfolding of $P$) and a morphism $f$ from $F$ onto $P$ which is a bijection between the minimal elements of $F$ and $P$, and for any non-minimal element $x\in F$, $f$ is a  bijection between the successors of $x$ in $F$ and the successors of $f(x)$ in $P$. Moreover, $F$ is a largest element in $(\{G\in{\mathcal F}_Q\mid
G\leq_hP\};\leq_h)$. This   extends to the iterated version: for any $(P,\leq,c)\in{\mathcal P}_{{\mathcal P}_Q}$ there is an $\leq_h$-largest element
$(F(P)\leq,d)\in{\mathcal F}_{{\mathcal F}_Q}$ below $P$:  it suffices to set 
$d=F\circ c\circ f$. If $c(p)$ has a least element for every $p\in P$ then $(F(P)\leq,d)\in{\mathcal F}_{{\mathcal T}_Q}$.

By a {\em minimal $Q$-forest} we mean a finite $Q$-forest not $h$-equivalent  to a $Q$-forest of lesser cardinality; this notion  also naturally extends to the iterated trees and forests. Proposition 8.3 in \cite{s12} gives an algorithm for checking minimality that also extends to the iterated labeled forests. The minimal $Q$-forests are useful because properties of elements of the quotient-poset of $(\mathcal{F}^\sqcup_k(\omega);\leq_h)$ are easier to check for the minimal representatives. In particular, the equivalence class $[F]_h$ is join-irreducible iff the minimal forest $M\equiv_hF$ is a tree.

For a further reference, we cite Proposition 8.13 from \cite{s12} that associates inductively (and effectively) to any $F\in\mathcal{T}_k(\omega)^\sqcup$, a  finite subset  $M(F)$ of $\mathcal{T}_k(\omega)$ such that $M(F)$ contains all (up to $\equiv_h$) minimal elements of $(\{G\in\mathcal{T}_k(\omega)^\sqcup\mid G\not\leq_h F\};\leq_h)$.

 \begin{Proposition}\label{min11q}
 \begin{enumerate}\itemsep-1mm
 \item For any $i<k$, $M(i)=\{j<k\mid  j\not=i\}$.
 \item If $F=F_0\sqcup\cdots\sqcup F_n$ is a minimal forest canonically decomposed to $n>0$ trees
 then $M(F)=\{j\cdot(G_0\sqcup\cdots\sqcup G_n)\mid j<k, G_0 \in M(F_0),\ldots,G_n \in M(F_n)\}$.
 \item If $F=i\cdot G$ is a minimal (in the sense explained above)  $k$-tree and $G$ is a nonempty forest then $M(F)=\{j\cdot K\mid j\in k\setminus\{i\},K\in M(G)\}$.
 \item If $F=s(V)$ and $V$ is a minimal tree in $\mathcal{T}_k(\omega)\setminus\overline{k}$ then
$M(s(V))=\{s(V_1)\mid V_1\in M(V)\}$.
 \item If $F=s(V)\cdot G$ is  minimal and $G$ is a nonempty forest
 then $M(F)=\{s(V_1)\cdot G_1\mid s(V_1)\in M(T),G_1\in M(G)\}$.
 \end{enumerate}
  \end{Proposition}

Next we recall another form of iterated labeled posets which  appear naturally in some contexts (we will see examples in the next section).
By {\em $\omega$-preorder} we mean a structure $C=(C;\leq_0,\leq_1,\ldots)$
where $\leq_n$ are preorders on $C$ such that
$x\leq_{n+1}y\rightarrow x\equiv_ny$, for each $n<\omega$. In this paper we mainly deal with the particular case of {\em
2-preorder} that is a structure $C=(C;\leq_0,\leq_1)$ with two
preorders on $P$ such that $x\leq_1y\rightarrow x\equiv_0y$ but some facts below are formulated for $n$-preorders $(C;\leq_0,\ldots,\leq_{n-1})$, $n>0$, defined in the same way.

A {\em $Q$-labeled 2-preorder} is a tuple
$(C;\leq_0,\leq_1,A)$ where $(C;\leq_0,\leq_1)$ is a 2-preorder and $A:C\to Q$. A {\em morphism of
$k$-labeled 2-preorders} $C,C_1$ is a function $g:C\to C_1$ that respects the preorders and satisfies $\forall x\in C(A(x)\leq A_1(g(x)))$. By Proposition 8.8 in  \cite{s12}, the category
 of finite $Q$-labeled $2$-preorders is equivalent to the category ${\mathcal P}_{{\mathcal P}_Q}$. The 2-iterated $Q$-poset corresponding to $(C;\leq_0,\leq_1,A)$ is $(C/_{\equiv_0};\leq_0,d)$  where 
$d([c]_0)=([c]_0;\leq_1,A|_{[c]_0})$ for every $c\in C$. In fact, this  holds for $Q$-labeled $n$-preorders for every $n>0$.

We conclude this subsection with answering the natural question about the complexity of first order theory $FO(\mathbb{F}_{k}(n))$ where $\mathbb{F}_{k}(n)$ is the quotient-poset of $(\mathcal{F}_{k}(n);\leq_h)$. The next theorem is an easy corollary of the result in \cite{ks07} that, for $k>2$, $FO(\mathbb{F}_{k}(1))$ is computably isomorphic to the first-order arithmetic. It shows that the poset $\mathbb{F}_{k}(n)$ for $k>2$ is really (not only intuitively from Figs. 1 and 2) much more complicated than for $k=2$.

\begin{Theorem}\label{fotheory}
Let $1\leq n\leq\omega$. If $k=2$ then $FO(\mathbb{F}_{k}(n))$ is decidable, otherwise it is undecidable and, moreover, computably isomorphic to the first-order arithmetic.  
 \end{Theorem}

{\em Proof.} The decidability for $k=2$ follows from the fact that the partial order $\mathbb{F}_{2}(n)$ is semi-well-ordered and the self-dual and non-self-dual levels alternate, so let $k>2$. For $n=1$ this is the result in \cite{ks07}, so let $n>1$. 

The result is deduced from the case $n=1$ and the following  facts contained in my papers on the iterated $h$-preorders (which  also easily follow from the remarks above):   $\mathbb{F}_{k}(1)$ is an initial segment of  $\mathbb{F}_{k}(n)$ consisting of all elements of finite rank, and the set of minimal elements in $\mathbb{F}_{k}(n)\setminus\mathbb{F}_{k}(1)$ contains precisely $k(k-1)$ elements. This leads to the following first-order definition of $\mathbb{F}_{k}(1)$ in $\mathbb{F}_{k}(n)$ without parameters (which is obviously sufficient for the proof).

Let $\lambda(u)$ be the formula $\exists y(y<u)\wedge\forall y<u\exists z(y<z<u)$ saying that $u$ is not
minimal and has no immediate predecessor. Let $\mu(u)$ be the formula $\lambda(u)\wedge\forall v<u\neg\lambda(u)$ saying that $u$ is minimal among the non-minimal elements having no immediate
predecessor. Finally, let $p=k(k-1)$,  $u_1,\ldots,u_p$ be different
variables, and $\phi(x)$ be the formula
$$\exists u_1\cdots\exists
u_p((\bigwedge_{i\not=j}u_i\not=u_j)\wedge(\bigwedge_i(\mu(u_i)\wedge
u_i\not\leq x))).$$
By the remarks in the previous paragraph, $\phi(x)$ defines
$\mathbb{F}_{k}(1)$ in $\mathbb{F}_{k}(n)$.
 \qed

\subsection{Bases and fine hierarchies}\label{fhpart}

Here we  briefly recall  technical notions  (slightly adapted from \cite{s12,s21}) related to the FHs.

By an {\em $\omega$-base in a set $S$} we mean a sequence $\mathcal{L}=\mathcal{L}(S)=\{\mathcal{L}_n\}_{n<\omega}$ of subalgebras of $(P(S);\cup,\cap,\emptyset,S)$ such that  $\mathcal{L}_n\cup\check{\mathcal{L}}_n\subseteq\mathcal{L}_{n+1}$ for each $n$. The $\omega$-base ${\mathcal L}$ is {\em reducible} if every ${\mathcal L}_n$ has the
reduction property (i.e., for
every $C_0, C_1\in{\mathcal L}$ there are disjoint
$C^\prime_0,C^\prime_1\in{\mathcal L}$ such that
$C^\prime_i\subseteq C_i$ for $i=0,1$, and $C_0\cup
C_1=C^\prime_0\cup C^\prime_1$).  The $\omega$-base ${\mathcal L}$ is {\em separable} if every $\check{\mathcal L}_{n+1}$ has the
separation property (i.e., every disjoint $\check{\mathcal L}_{n+1}$-sets $A,B$ are separable by a suitable $C\in{\mathcal L}_{n+1}\cap\check{\mathcal L}_{n+1}$: $A\subseteq C\subseteq\overline{B}$). 
The $\omega$-base ${\mathcal L}$ is {\em interpolable} if for every $n<\omega$ any two
disjoint sets in $\check{\mathcal L}_{n+1}$ are separable by a set in $BC({\mathcal L}_n)$. Clearly, every reducible base is separable, and the base ${\mathcal L}$ is interpolable iff it is separable and ${\mathcal L}_{n+1}\cap\check{\mathcal L}_{n+1}=BC({\mathcal L}_n)$ for every $n<\omega$. 

Similarly to the previous subsection, in this paper we mainly deal with particular cases of $\omega$-bases, namely with 1- and 2-bases. A {\em 1-base} in  $S$ is just a subalgebra of $\mathcal{L}$ of $(P(S);\cup,\cap,\emptyset,S)$. A {\em 2-base} in  $S$ is a pair $\mathcal{L}=(\mathcal{L}_0,\mathcal{L}_1)$ of 1-bases in $X$ such that $\mathcal{L}_0\subseteq\mathcal{L}_1\cap\check{\mathcal{L}}_1$.

With any $\omega$-base $\mathcal{L}(S)$  we  associate some other $\omega$-bases as follows. For any $m<\omega$, let $\mathcal{L}^m(S)=\{\mathcal{L}_{m+n}(S)\}_n$; we call this base  the {\em $m$-shift of $\mathcal{L}(S)$}.  For any $U\in\mathcal{L}_0$, let $\mathcal{L}(U)=\{\mathcal{L}_n(U)\}_{n<\omega}$ where $\mathcal{L}_n(U)=\{U\cap V\mid V\in\mathcal{L}_{n}(S)\}$; we call this base in $U$ the {\em $U$-restriction of $\mathcal{L}(S)$}. For any  subalgebra $\mathcal{B}$ of $(P(S);\cup,\cap,\bar{\;},\emptyset,S)$, let $\mathcal{B}\mathcal{L}(S)=\{\mathcal{B}\cap\mathcal{L}_{n}(S)\}_n$; we call this base in $S$ the {\em $\mathcal{B}$-fragment of $\mathcal{L}(S)$}. 

With any $\omega$-base $\mathcal{L}(S)$ in $S$ we associate the {\em FH of $k$-partitions over $\mathcal{L}$} that is a family $\{\mathcal{L}(S,T)\}_{T\in\mathcal{T}_k(\omega)}$ of subsets of $k^S$. Note that, for $k=2$, we obtain the FH of sets over $\mathcal{L}$ that, according to the structure of $(\mathcal{T}_2(\omega);\leq_h)$ in the previous subsection, looks much easier: it may be written as the sequence 
$\{\mathcal{L}(\alpha)\}_{\alpha<\varepsilon_0}$ whose members correspond to the $\Sigma$-levels. In fact, we will define $\mathcal{L}(S,T)$ not only for $T\in\mathcal{T}_k(\omega)$ but also for $T\in\mathcal{T}_k(\omega)^\sqcup$.

To avoid some technicalities, we give definitions here only for the 1- and 2-bases; for the general case see e.g. Section 3 of \cite{s21}. Let first $\mathcal{L}(S)$ be a 1-base and $(F,c)\in\mathcal{F}_k$, $F\subseteq\omega^+$.
 An {\em $F$-family over $\mathcal{L}$} is a family $\{U_\tau\}_{\tau\in F}$ of $\mathcal{L}$-sets such that  $U_\tau\supseteq U_{\tau i}$ for $\tau i\in F$ and $\bigcup_{\tau\in F}U_\tau=S$. If, in addition, $U_{\tau i}\cap U_{\tau j}=\emptyset$ for all $\tau i,\tau j\in F$ with $i\neq j$, the family $\{U_\tau\}_{\tau\in F}$ is called {\em reduced}. We say that  the $F$-family {\em determines a $k$-partition $A:S\to\bar{k}$} if $A(x)=c(\tau)$ where $\tau$ is (any) string in $F$ such that $x\in\tilde{U}_\tau=U_\tau\setminus\bigcup\{U_{\tau i}\mid\tau i\in F\}$. Note that any $F$-family determines at most one $k$-partition, and any reduced $F$-family determines precisely one $k$-partition (because in this case the {\em components} $\tilde{U}_\tau$ form a partition of $S$). The {\em FH of $k$-partitions over $\mathcal{L}$} is the family $\{\mathcal{L}(T)\}_{T\in\mathcal{F}_k}$ where $\mathcal{L}(T)$ is the set of $k$-partitions determined by $T$-families over $\mathcal{L}$.  

Let now $\mathcal{L}$ be a 2-base and $(F,c)\in\mathcal{T}_k(2)^\sqcup$, $F\subseteq\omega^+$. An {\em $F$-family over $\mathcal{L}$} is a pair $(\{U_\tau\},\{U_{\tau\sigma}\})$ where $\{U_\tau\}$ is an $F$-family over $\mathcal{L}_0$ and, for any $\tau\in F$, $\{U_{\tau\sigma}\}$ is a $c(\tau)$-family over the shifted 1-base $\{\widetilde{U}_\tau\cap B\mid B\in\mathcal{L}_1\}$ in $\widetilde{U}_\tau$. The family $(\{U_\tau\},\{U_{\tau\sigma}\})$ is {\em reduced} if  $\{U_\tau\}$ and $\{U_{\tau\sigma}\}$ for each $\tau$ are reduced. The $F$-family over $\mathcal{L}$ {\em determines a $k$-partition $A:S\to\bar{k}$} if $A(x)=v(\sigma)$, for every $\tau\in F$ and $\sigma\in c(\tau)=(V,v)\in\mathcal{T}_k$ such that $x\in\widetilde{U}_{\tau\sigma}$. Again, any family determines at most one $k$-partition, and any reduced family determines precisely one $k$-partition (because in this case the {\em final components} $\widetilde{U}_{\tau\sigma}$ form a partition of $X$). The {\em FH of $k$-partitions over the 2-base $\mathcal{L}$} is the family $\{\mathcal{L}(F)\}_{F\in\mathcal{F}_{\mathcal{T}_k}}$ where $\mathcal{L}(F)$ is the set of $k$-partitions determined by $F$-families over $\mathcal{L}$.

We summarize some properties of the FH $\{\mathcal{L}(S,F)\}$ over an $\omega$-base $\mathcal{L}$ in $S$. 

\begin{Proposition}\label{fhkp}
\begin{enumerate}\itemsep-1mm
\item If $F\leq_hG$ then $\mathcal{L}(S,F)\subseteq\mathcal{L}(S,G)$ (and hence $(\{\mathcal{L}(S,F)\mid F\in\mathcal{T}_k(\omega)^\sqcup\};\subseteq)$ is a wqo).
\item If $\mathcal{L}$ is reducible then every level $\mathcal{L}(S,F)$ coincides with the set of $k$-partitions determined by the reducible $F$-families over $\mathcal{L}$.
\item If $\mathcal{L}$ is interpolable then the FH of sets over $\mathcal{L}$ has the following discreteness property: for every limit ordinal $\lambda<\varepsilon_0$, $\mathcal{L}(\lambda)\cap\check{\mathcal{L}}(\lambda)=\bigcup\{\mathcal{L}(\alpha)\mid\alpha<\lambda\}$.
\item If $\mathcal{L}'$ is another $\omega$-base in $S'$ and $f:S\to S'$ satisfies $\forall n\forall A\in\mathcal{L}'_n(f^{-1}(A)\in\mathcal{L}_n)$ then $A\circ f\in\mathcal{L}(F)$ for all $A\in\mathcal{L}'(F)$ and $F\in\mathcal{T}_k(\omega)^\sqcup$.
 \item For every tree $T$, $\mathcal{L}(S,s(T))=\mathcal{L}^1(S,T)$ 
 \item Let $F=F_0\sqcup\cdots\sqcup F_n$, $n>0$, be a minimal proper forest with trees $F_i$, let $U_0,\cdots,U_n$ be $\mathcal{L}_0$-sets that cover $S$, and let $A\in k^S$. If $A|_{U_i}\in\mathcal{L}(U_i,F_i)$ for all $i\leq n$ then $A\in\mathcal{L}(S,F)$.
 \item Let $s(T)\cdot F$ be a minimal forest with non-empty $F$,  $A\in k^S$, and let $U\in\mathcal{L}_0$ satisfy $A|_{\overline{U}}\in\mathcal{L}^1(\overline{U},T)$ and $A|_U\in\mathcal{L}(U,F)$. Then $A\in\mathcal{L}(S,s(T)\cdot F)$.
\item For every subalgebra $\mathcal{B}$ of $(P(S);\cup,\cap,\bar{\;},\emptyset,S)$ and for every $F\in\mathcal{F}_k(\omega)$, we have: $\mathcal{B}\mathcal{L}(F)\subseteq\mathcal{B}_k\cap\mathcal{L}(F)$.
 \end{enumerate}
 \end{Proposition}
 
 {\em Proof Hint.} The proofs of (1)---(7) may be found in  \cite{s12,s21}, so we check only item (8) that is formally new. Let $A\in\mathcal{B}\mathcal{L}(F)$, then $A$ is determined by an $F$-family $(\{U_\tau\},\{U_{\tau\sigma}\})$ over $\mathcal{B}\mathcal{L}=(\mathcal{B}\cap\mathcal{L}_0,\mathcal{B}\cap\mathcal{L}_1)$ (to simplify notation, we take $F\in\mathcal{F}_{\mathcal{T}_k}$), so in particular, $U_\tau\in\mathcal{B}\cap\mathcal{L}_0$ and $U_{\tau\sigma}\in\mathcal{B}\cap\mathcal{L}_1$. Then $(\{U_\tau\},\{U_{\tau\sigma}\})$ is also an $F$-family over $\mathcal{L}$, hence $A\in\mathcal{L}(F)$. Since every $\tilde{U}_\tau$ (resp. $\tilde{U}_{\tau\sigma}$) is a difference of $\mathcal{B}\cap\mathcal{L}_0$-sets (resp. of $\mathcal{B}\cap\mathcal{L}_1$-sets), every component $A^{-1}(i)$, $i<k$, is in $\mathcal{B}$. Therefore, $A\in\mathcal{B}_k$.
 \qed

Item (1) motivates the following definition: the FH over $\mathcal{L}$ {\em does not collapse} if $\mathcal{F}_k(\omega)\simeq(\{\mathcal{L}(X,F)\mid F\in\mathcal{T}_k(\omega)^\sqcup\};\subseteq)$. Item (8) motivates the following definition: the Boolean algebra $\mathcal{B}$ is {\em $\mathcal{L}$-precise}, if $\mathcal{B}\mathcal{L}(F)=\mathcal{B}_k\cap\mathcal{L}(F)$ for every $F\in\mathcal{F}_k(\omega)$. In general, precise Boolean algebras are rare but below we will show that both $\mathcal{R}$ and $\mathcal{A}$ are precise w.r.t. the Borel hierarchy.

\subsection{Examples of fine hierarchies}\label{sample}

Here we discuss examples of FHs most relevant to this paper.

The first example is the {\em  Borel base} $\{\bfSig^0_{1+n}\}$ in  $X^\omega$. This $\omega$-base is  known to be reducible but it is not interpolable because, for each $n<\omega$,  $\bfDe^0_{2+n}$ is much larger than $BC(\bfSig^0_{1+n})$. The FH of $k$-partitions over this base could be denoted by $\{\bfSig(F)\}_{F\in\mathcal{F}_k(\omega)}$. The FH of sets over this base (for the non-self-dual levels) could be denoted as $\{\bfSig_\alpha\}_{\alpha<\varepsilon_0}$.
Since the class $\mathcal{R}$ of regular sets in $X^\omega$ is contained in $BC(\bfSig^0_{2})$, the 2-base $\bfSig=(\bfSig^0_1,\bfSig^0_2)$ is in fact sufficient for the sequel. The FH of $k$-partitions  over $\bfSig$ looks as $\{\bfSig(F)\}_{F\in\mathcal{F}_k(2)}$. The FH of sets over this base (for the non-self-dual levels) could be written as $\{\bfSig_\alpha\}_{\alpha<\omega^\omega}$.

The second example is the 2-base $\mathcal{R}\bfSig$ that is the $\mathcal{R}$-fragment of $\bfSig$. As shown in \cite{s98}, the 2-base $\mathcal{R}\bfSig$ is reducible. It is also interpolable because the interpolability of this base is equivalent to the Staiger-Wagner theorem \cite{sw74}.  The FH of $k$-partitions  over this base will be denoted by $\{\mathcal{R}\bfSig(F)\}_{F\in\mathcal{F}_k(2)}$. The FH of sets over this base (for the non-self-dual levels)  was denoted as $\{\mathcal{R}_\alpha\}_{\alpha<\omega^\omega}$ in \cite{s98}.

The third example is the 2-base $\mathcal{A}\bfSig$ that is the $\mathcal{A}$-fragment of $\bfSig$, where $\mathcal{A}$ is the class of regular aperiodic sets in $X^\omega$. As shown in \cite{s08}, the 2-base $\mathcal{A}\bfSig$ is reducible and interpolable. The FH over of $k$-partitions  this base will be denoted by $\{\mathcal{A}\bfSig(F)\}_{F\in\mathcal{F}_k(2)}$. The FH of sets over this base (for the non-self-dual levels)  was denoted as $\{\mathcal{A}_\alpha\}_{\alpha<\omega^\omega}$ in \cite{s08}.

The forth example is the base $\{\mathcal{C}_n\}_{n<\omega}$ associated with any $\omega$-preorder $(C;\leq_0,\ldots)$ where, for each $n<\omega$, $\mathcal{C}_n$ is the class of all $\leq_n$-up subsets of $C$ (a set $A\subseteq C$ is $\leq_n$-up if $a\in A$ and $a\leq_nc$ imply $c\in A$). Such bases and the corresponding FH's $\{\mathcal{C}(F)\}_{F\in\mathcal{F}_k(\omega)}$  appear in the sequel; more precisely, in this paper we mainly deal with 2-preorders $(C;\leq_0,\leq_1)$ and the corresponding FH's $\{\mathcal{C}(F)\}_{F\in\mathcal{F}_k(2)}$ over the 2-bases $\mathcal{C}=(\mathcal{C}_0,\mathcal{C}_1)$.
We formulate some properties of such bases and FH's depending on the properties of the source $\omega$- or 2-preorders.

 We call a 2-preorder $(C;\leq_0,\leq_1)$ {\em compatible} if $a\equiv_0b$ implies $\exists c(c\leq_1a\wedge c\leq_1b)$; note that the converse implication holds in every 2-preorder. For  a subset $A$ of a preorder $(C;\leq)$, let $\downarrow{A}=\{x\mid\exists a\in A(x\leq a)\}$ and $\uparrow{A}=\{x\mid\exists a\in A(a\leq x)\}$. If $A=\{a\}$ is a singleton, we simplify $\downarrow{A}$ to $\downarrow{a}$ and $\uparrow{A}$ to $\uparrow{a}$.

\begin{Lemma}\label{compat}
If the 2-preorder $(C;\leq_0,\leq_1)$ is compatible then $\check{\mathcal{C}}_1$ has the separation property. If, in addition, $C$ is finite then the 2-base $(\mathcal{C}_0,\mathcal{C}_1)$ is interpolable.
 \end{Lemma} 

{\em Proof.} Let $A,B$ be disjoint sets in $\check{\mathcal{C}}_1$. Let $[A]_0$ be the closure of $A$ under $\equiv_0$, then $A\subseteq[A]_0$ and $[A]_0$ is in $\mathcal{C}_1\cap\check{\mathcal{C}}_1$, hence for the first assertion it remains to show that $[A]_0$ is disjoint from $B$. Suppose the contrary: $b\in[A]_0$ for some $b\in B$. Then $a\equiv_0b$ for some $a\in A$ (the case $A=\emptyset$ is trivial). By compatibility,  $c\leq_1a,b$ for some $c\in C$. Since both of $A,B$ are $\leq_1$-down, $c\in A\cap B$. Contradiction.

Let now $C$ be finite. To prove the second assertion, it suffices to show that $[A]_0$ is in $BC(\mathcal{C}_0)$. Since $[A]_0$ is the finite union of the equivalence classes $[a]_0$ for all $a\in A$, it suffices to show that $[a]_0$ is in $BC(\mathcal{C}_0)$. Clearly, $[a]_0=\uparrow_0a\cap\downarrow_0a$ is a difference of $\mathcal{C}_0$-sets. 
 \qed

\begin{Corollary}\label{compat1}
If the $\omega$-preorder $(C;\leq_0,\ldots)$ is such that the 2-preorder $(C;\leq_n,\leq_{n+1})$ is compatible for every $n<\omega$ then the $\omega$-base  $\mathcal{C}$ is separable. If, in addition, $C$ is finite then   $\mathcal{C}$ is interpolable.
 \end{Corollary}
 
The next lemma is a preparation to the subsequent theorem. 

\begin{Lemma}\label{pre1}
Let $(C;\leq)$ be a  1-preorder and $A,B_0,\ldots,B_n\subseteq C$ satisfy
   $\forall a\in A(\forall b\in B_0(a\not\leq b)\vee\cdots\vee\forall b\in B_n(a\not\leq b)$.
Then there exist  $U_0,\ldots U_n\in\mathcal{C}$ such that
$A\subseteq U_0\cup\cdots\cup U_n$ and $A\cap
U_i\subseteq\overline{B}_i$ for all $i\leq n$.
\end{Lemma}

{\bf Proof.} Let $A_i=\{a\in A\mid\forall b\in B_i(a\not\leq b)\}$ for each $i\leq n$, so $A=A_0\cap\cdots\cup A_n$. Then it suffices to take $U_i=\uparrow{A}_i$.
 \qed

Though the following theorem is just a technical characterisation of the FH over $n$-preorders in terms of the $h$-preorder, it is quite useful in the sequel. It extends or resembles some  characterisations of difference hierarchies in terms of the so called alternating chains \cite{ad65,s08a}, later extended to characterisations of FHs of sets in terms of alternating trees \cite{s95,s08a}, and then to characterisations of  difference and FHs of $k$-partitions (cf. Theorems 7.18 and 8.23 in \cite{s12}).

\begin{Theorem}\label{alttree}
Let $(C;\leq_0,\ldots,\leq_{n-1})$ be an $n$-preorder, $0<n<\omega$, $A\in k^C$, and $F\in\mathcal{F}_k(n)$. Then $A\in\mathcal{C}(F)$ iff $T\leq_hF$ for every $T\in\mathcal{T}_k(n)$ with $T\leq_h (C;\leq_0,\ldots,\leq_{n-1},A)$. 
 \end{Theorem}

{\em Proof.} The proof is by induction on $n$ but for simplicity we consider here only the cases $n=1,2$ sufficient for this paper. The  direction from left to right will be checked for arbitrary preorder $Q$ in place of $\bar{k}$. For $(F,f)\in\mathcal{F}_Q$, we say that {\em $A$ is majorised by an $F$-family $\{U_\rho\}_{\rho\in F}$ over $\mathcal{C}$} if $A(c)\leq f(\rho)$ for every $\rho\in F$ with $c\in\tilde{U}_\rho$. For $F\in\mathcal{F}_{\mathcal{T}_Q}$, we say that $A$ is {\em  majorised by an $F$-family $(\{U_\rho\},\{U_{\rho\sigma}\}_{\sigma\in f(\rho)})$ over $\mathcal{C}$} if $A(c)\leq v(\sigma)$ for every $\rho\in F$ and $\sigma\in f(\rho)=(V,v)$ with $c\in\tilde{U}_{\rho\sigma}$.

For $n=1$, it suffices to show that if $A$ is majorised by an $F$-family $\{U_\tau\}$ and $T\leq_h(C;\leq_0,A)$ via a monotone function $\varphi:T\to(C;\leq_0)$ (where $T\in\mathcal{T}_Q$) then $T\leq_hF$. It suffices to construct by induction a monotone function $\psi:T\to F$ such that $t(\tau)\leq A(\varphi(\tau))\leq f(\psi(\tau))$ and $\varphi(\tau)\in\tilde{U}_{\psi(\tau)}$ for every $\tau\in T$. For some $\rho\in F$ we have $\varphi(\varepsilon)\in\tilde{U}_\rho$; let $\psi(\varepsilon)$ be any such $\rho$, then the desired condition holds. Assume by induction that it holds for a fixed $\tau$ and let $\tau i\in T$, so $t(\tau i)\leq A(\varphi(\tau i))$. Since $\varphi(\tau)\in U_{\psi(\tau)}$, $\varphi(\tau)\leq_0\varphi(\tau i)$, and $U_\tau$ is $\leq_0$-up, we get $\varphi(\tau i)\in U_{\psi(\tau)}$. Then $\varphi(\tau i)\in\tilde{U}_\rho$ for some $\rho\in F$ with $\rho\sqsupseteq\psi(\tau)$; setting $\psi(\tau i)=\rho$  for some such $\rho$ completes the induction step.

For $n=2$, it suffices to show that if $A$ is majorised by an $F$-family $(\{U_\rho\},\{U_{\rho\sigma}\})$ and $T\leq_h(C;\leq_0,\leq_1,A)$  (where $T\in\mathcal{T}_{\mathcal{T}_Q}$ and $F\in\mathcal{F}_{\mathcal{T}_Q}$) then $T\leq_hF$. Note that the notation here is slightly abused (instead of $(C;\leq_0,\leq_1,A)$ should stay its copy $(C/_{\equiv_0};\leq_0,d)$ from $\mathcal{P}_{\mathcal{P}_Q}$ where $d:C/_{\equiv_0}\to\mathcal{P}_Q$ is given by $d([c]_0)=([c]_0\leq_1,A|_{[c]_0})$. Thus, a witness for $T\leq_h(C;\leq_0,\leq_1,A)$ may be given by a monotone function $\varphi:(T\sqsubseteq)\to(C;\leq_0)$ such that $t(\tau)\leq_hd([\varphi(\tau)]_0)$ for each $\tau\in T$. For every $\tau\in T$ there are $\rho,\sigma$ with $\varphi(\tau)\in\tilde{U}_{\rho\sigma}$, and for every such $\rho,\sigma$ we have $A(\varphi(\tau))\leq v(c)$. Since $\varphi(\tau)\in\tilde{U}_{\rho}$ and $U_\rho$ is an $\mathcal{L}_0$, $[\varphi(\tau)]_0\subseteq\tilde{U}_{\rho}$. Since $A|_{[\varphi(\tau)]_0}$ is majorised by the $f(\rho)$-family $\{[\varphi(\tau)]_0\cap U_{\rho\sigma}\}$ over the restricted base $\mathcal{C}_1([\varphi(\varepsilon)]_0)$ and $t(\tau)\leq_h([\varphi(\tau)]_0\leq_1,A|_{[\varphi(\tau)]_0})$, by the case $n=1$ (applied to the restricted 1-base $\mathcal{C}_1([\varphi(\tau)]_0)$ in place of $\mathcal{L}_0$, $t(\tau)$ in place of $T$, $f(\rho)$ in place of $F$, and $A|_{[\varphi(\tau)]_0}$ in place of $A$) we get $t(\tau)\leq_hf(\rho)$. Now we can repeat the induction over $T$ from the previous paragraph and construct a monotone function $\psi:T\to F$ such that $t(\tau)\leq_hf(\psi(\tau))$ for every $\tau\in T$. Thus, $T\leq F$.  

In the other direction, let $n=2$ and $T\leq_hF$ for all $T\in\mathcal{T}_{\mathcal{T}_k}$ with $T\leq_h(C;\leq_0,\leq_1,A)$. This implies that $T\not\leq_h(C;\leq_0,\leq_1,A)$ whenever $T\in M(F)$. For any $G\in\mathcal{F}_{\mathcal{T}_k}$, let $B_G$ be the set of all $c\in C$ such that there is a morphism $\varphi:G\to (C;\leq_0,\leq_1,A)$ with $c\leq_0\varphi(G)$ (i.e., $c\leq_0d$ for every $d\in\varphi(G)$); note that $B_G\in\check{\mathcal{C}_0}$. For any $G\in\mathcal{F}_{k}$, let $B^1_G$ be the set of all $c\in C$ such that there is a morphism $\varphi:G\to (C;\leq_1,A)$ with $c\leq_1\varphi(G)$; note that $B^1_G\in\check{\mathcal{C}_1}$.

We prove $A\in\mathcal{L}(F)$ by induction on the cardinality of the forest $F$ (assuming w.l.o.g. that $F$ is minimal in the sense of Subsection \ref{lab}). For singleton forest $F=i<k$ the assertion follows immediately from  Proposition \ref{min11q}(1). For the non-trivial singleton tree $F=s(V)$, the assertion reduces to the case $n=1$ by Propositions \ref{min11q}(4) and \ref{fhkp}(5).

Let now  $|F|\geq2$ and $F$ be not a tree, i.e.
$F=F_0\sqcup\cdots\sqcup F_m$ for some $m\geq1$ and minimal
$k$-trees $F_0,\ldots,F_n$. By  Proposition \ref{min11q}(2),
for all $G_0\in M(F_0),\ldots,G_m\in M(F_m)$ and $j< k$ we have
$j\cdot(G_0\sqcup\cdots\sqcup G_m)\not\leq_h(C;\leq_0,\leq_1,A)$. For
any $l\leq m$, let $B_l=\bigcup\{B_G\mid G\in M(F_l)\}$. Then
$\forall a\in C(\forall b\in B_0(a\not\leq_0 b)\vee\cdots\vee\forall
b\in B_m(a\not\leq_0 b))$. Suppose the contrary: there are  $a\in C$
and $b_0\in B_0,\ldots,b_m\in B_m$ with $a\leq_0 b_0,\ldots,b_m$.
For any $l\leq m$, choose $G_l\in M(F_l)$ with $b_l\in B_{G_l}$,
then $a\leq_0\varphi_l(G_l)$ for some morphism
$\varphi_l:G_l\to(C;\leq_0\leq_1,A)$. Define a mapping
$\varphi:H\rightarrow C$ (where $H=j\cdot(G_0\sqcup\cdots\sqcup
G_m)$  and $j= A(a)$) by: $\varphi(\varepsilon)=a$  and $\varphi(l\tau)=\varphi_l(\tau)$ for all $l\leq
m$ and $\tau\in G_l$. Then $\varphi$ is a morphism from $H$
to $(C;\leq_0,\leq_1,A)$, which is a contradiction. By Lemma
\ref{pre1}, there are $\mathcal{C}_0$-sets $U_0,\ldots,U_m$ such
that $C\subseteq U_0\cup\cdots\cup U_m$ and
$U_l\subseteq\overline{B}_l$ for each $l\leq n$.  Then for the
$k$-partitions $A|_{U_l}$, $l\leq m$, we have $\forall
K\in M(F_l)(K\not\leq(U_l;\leq_0,\leq_1,A|_{U_l})$ because $U_l$
are $\leq_0$-up-sets (otherwise, $\varphi(\varepsilon)\in U_l\cap B_l$ for the witnessing morphism $\varphi$). Thus, by the definition of $M(F_l)$, the right-hand side condition in the formulation of the theorem hold for $F_l$ in place of $T$, and $(U_l;\leq_0,\leq_1,A|_{U_l})$ in place of $(C;\leq_0,\leq_1,A)$. By induction,
$A|_{U_l}\in\mathcal{L}(U_l,F_l)$ for each $l\leq m$. By 
Proposition \ref{fhkp}(6), $A\in\mathcal{C}(F)$.

Finally, let $F=s(V)\cdot G$ be a minimal tree where $V\not\in\overline{k}$ and
$G\not=\emptyset$ (the case $V=i<k$ is considered similarly, using Proposition \ref{min11q}(3)). Since, by Proposition \ref{min11q}(5),
$s(V_1)\cdot G_1\not\leq_h s(V)\cdot G$ for all $V_1\in M(V),G_1\in M(G)$, we have
$(s(V_1)\cdot G_1)\not\leq_h(C;\leq_0,\leq_1,A)$ for
all such $V_1,G_1$. The $\check{\mathcal{C}_1}$-set $B^1=\bigcup\{B^1_{s(V_1)}\mid V_1\in
M(V)\}$ and the $\check{\mathcal{C}_0}$-set $D=\bigcup\{B_{G_1}\mid G_1\in M(G)\}$ satisfy $\forall b\in
B\forall d\in D(b\not\leq_0 d)$. Suppose the contrary: $b\leq_0 d$
for some $b\in B,d\in D$. Then $b\leq_1\psi(V_1)$ for some
$V_1\in M(V),\psi:V_1\rightarrow(C;\leq_1,A)$, and $d\leq_0\theta(G_1)$ for
some $G_1\in
M(G),\theta:G_1\rightarrow(C;\leq_0,\leq_1,A)$,
and one easily, similar to the previous paragraph, constructs a morphism $\varphi:(s(V_1)\cdot
G_1)\rightarrow(C;\leq_0,\leq_1,A)$, which is a
contradiction. By Lemma \ref{pre1} (taken for $n=1$), $B\subseteq
U\subseteq\overline{D}$ for some $U\in\mathcal{C}_0$. Then
$\forall V_1\in
M(V)(V_1\not\leq_h(\overline{U};\leq_1,A|_{\bar{U}}))$ and $\forall G_1\in
M(G)(G_1\not\leq_h(U;\leq_0,\leq_1, A|_U))$ because
$U\in\mathcal{C}_0\subseteq\mathcal{C}_1\cap\check{\mathcal{C}}_1$.
By induction, $A|_{\bar{U}}\in\mathcal{C}^1(\bar{U},P)$ and
$A|_U\in\mathcal{C}(U,Q)$. By   Proposition \ref{fhkp}(7),
$A\in\mathcal{L}(F)$.
 \qed

\section{Classifying regular (aperiodic) $k$-partitions}\label{regular}

In this section we prove  main results of this paper. A basic fact is the following characterisation of Wadge degrees of regular (aperiodic) $k$-partitions.

\begin{Theorem}\label{wadred}
The quotient-posets of $(\mathcal{R}_k;\leq_{CA})$, $(\mathcal{A}_k;\leq_{CA})$, and $({\mathcal F}_{{\mathcal T}_k};\leq_h)$ are isomorphic.
 \end{Theorem}

The proof  below is in a sense an ``automatic'' constructive version of the corresponding proof in Wadge theory \cite{s17,km19}. It is divided into separate parts.

\subsection{Embedding ${\mathcal F}_{{\mathcal T}_k}$ into $\mathcal{A}_k$}\label{reduc0}

Here we define a function $\rho:{\mathcal F}_{{\mathcal T}_k}\to\mathcal{A}_k$ that will be shown to induce the desired isomorphisms in Theorem \ref{wadred}.

We start with introducing some operations on $k$-partitions.
The binary operation $A\oplus B$ on $k^{X^\omega}$ is defined by: $(A\oplus B)(0\xi)=A(\xi)$ and $(A\oplus
B)(i\xi)=B(\xi)$ for all $0<i<m$ and $\xi\in X^\omega$ (recall
that $X=\{0,\ldots,m-1$). Then $(k^{X^\omega};\leq_{CA},\oplus)$ is a $d$-semilattice. Since ${\mathcal R}_k$ and ${\mathcal A}_k$ are closed under $\oplus$, the structures $(\mathcal{A}_k;\leq_{CA},\oplus)$ and $(\mathcal{R}_k;\leq_{AS},\oplus)$ are semilattices (from Theorem \ref{wadred} it will follow that they are in fact also $d$-semilattices). 

The remaining facts are about some operations on $k^{X^\omega}$  closely related to the corresponding operations in \cite{s11,s17};  modifications are designed to make the set of aperiodic $k$-partitions closed under these operations. First,  we recall the unary operations $q_0,\ldots,q_{k-1}$ on
$k^{X^\omega}$ from \cite{s11} (that extend and modify the operation $\#$ from
\cite{wad84}) which use a coding of alphabets to guarantee the preservation of aperiodicity.  To simplify notation, we do this only for the
binary alphabet $X=\{0,1\}=2$ (however, it will be clear how to modify the idea for larger alphabets). Define the function $f:3^\omega\to 2^\omega$ by $f(x_0x_1\cdots)=\tilde{x}_0\tilde{x}_1\cdots$ where
$x_0,x_1\ldots<3$ and $\tilde{0}=110000,\tilde{1}=110100,\tilde{2}=110010$ (in the same
way we may define $f:3^\ast\to 2^\ast$). It is easy to see that $f$ is an $AA$-function, its image
$f(3^\omega)$ is a closed aperiodic set, and there is an
$AA$-function $f_1:2^\omega\to 3^\omega$ such that $f_1\circ
f=id_{3^\omega}$. For all $i<k$ and  $A\in k^{X^\omega}$, we define $q_i(A)\in k^{X^\omega}$ by
\[
[q_i(A)](\xi)=\left \{ \begin{array}{lll} i, & \mbox{if} &
\xi\not\in f(3^\omega)\vee\forall p\exists n\geq p(\xi[n,n+5]=\tilde{2}),\\
A(f_1(\xi)),& \mbox{if}& \xi\in f(2^\omega),\\
A(\eta),& \mbox{if}& \xi=f(\sigma2\eta)
\end{array}
\right.
\]
for some (unique) $\sigma\in 3^\omega$ and $\eta\in 2^\omega$. 

Using the same coding of $3^\omega$ into $2^\omega$, we define the binary operation  $\cdot$ on $k^{X^\omega}$ (that is a modification of the operation $+$ from \cite{wad84,s11}) as
follows (we again consider the typical particular case $X=2$). Define an $AA$-function $g:X^\omega\to X^\omega$ by
$g(x_0x_1\cdots)=\tilde{x}_0\tilde{2}\tilde{x}_1\tilde{2}\cdots$ where $x_0,x_1,\ldots\in X$ (in the same
way we may define  $g:X^\ast\to X^\ast$). Obviously,
$g(X^\omega)$ is a closed aperiodic set and there is an
$AA$-function $g_1:X^\omega\to X^\omega$ such that $g_1\circ
g=id_{X^\omega}$. For all $k$-partitions $A,B$ of  $X^\omega$, we
set
\[
[A\cdot B](\xi)=\left \{ \begin{array}{lll} A(g_1(\xi)), & \mbox{if} &
\xi\in g(X^\omega),\\
B(\eta),& \mbox{if} & \xi=g(u)\cdot v\cdot\eta,
\end{array}
\right.
\]
where $u\in X^*,\eta\in
X^\omega$, and $v\in X^+$ is the shortest word such that $g(u)\cdot v\cdot X^\omega\cap g(X^\omega)=\emptyset$.  
We will also use unary operations $p_i(A)\equiv_{AA}\bi\cdot A$ where $\bi=\lambda x.i$, $i<k$, is the constant $k$-partition. Equivalent operations with the names $p_0,\ldots,p_{k-1}$ were used  in \cite{s11,s17}.

\begin{Lemma}\label{closed}
The classes $\mathcal{A}_k$ and  $\mathcal{R}_k$ are closed under the  operations $q_0,\ldots,q_{k-1},\cdot$.
\end{Lemma}

{\em Proof.} We consider only  $\mathcal{A}_k$ but the proof also works for  $\mathcal{R}_k$. 
 Let $B=q_i(A)$ and let $A$ be aperiodic. It suffices to prove that $B_j$ is aperiodic for any $j\in\bar{k}\setminus\{i\}$. By the definition, $\xi\in B_j$ iff $\xi\in f(3^\omega)$ and there are only finitely many $n$ with $\xi[n,n+6)=\tilde{2}$ and either ($\xi\in f(2^\omega)$ and $f_1(\xi)\in A_j$) or $\exists n(\xi[n,n+6)=\tilde{2}\wedge\forall m>n(\xi[n,n+6)\not=\tilde{2})\wedge\xi[n+6,\infty)\in f_1^{-1}(A_j))$.  Since $A_j$ is aperiodic, so is also $B_j$, by the logical characterisation of aperiodic sets.

For the operation $\cdot$, let $C=A\cdot B$ and $i<k$. Then $\xi\in C_i$ iff ($\xi\in g(X^\omega)$ and $g_1(\xi)\in A_i$) or  ($\xi\not\in g(X^\omega)$ and $\eta\in B_i$) where $\eta$ as in the definition of $C$. From aperiodicity of $A_i,B_i$ and the definition of $u,v$ it is easy to find a first-order sentence defining $C_i$. Thus, $C$ is aperiodic.
 \qed

Next we define functions $\mu,\nu,\rho$ (which are variants of the corresponding functions from Section 6 of \cite{s17}) from labeled trees to $k$-partitions using the operations $p_i,q_i,\cdot$. 
Let $(T;t)\in\mathcal{T}_k$ where $T\subseteq\omega^*$. We associate with any  $\tau\in T$ the $k$-partition $\mu_T(\tau)$ by induction on the rank $rk(\tau)$ of $\sigma$ in $(T;\sqsupseteq)$ as follows: if $rk(\tau)=0$, i.e. $\tau$ is a leaf of $T$ then $\mu_T(\tau)=\bf{i}$ where $i=t(\tau)$; otherwise, $\mu_T(\tau)=p_i(\bigoplus\{\mu_T(\tau n)\mid n<\omega,\tau n\in T\})$. 
Now we define a function $\mu:\mathcal{T}_k\to k^{X^\omega}$ by $\mu(T)=\mu_T(\varepsilon)$. We define $\nu:\mathcal{T}_k\to k^{X^\omega}$ in the same way but using $q_i$ instead of $p_i$.

Now let $(T;t)\in\mathcal{T}_{\mathcal{T}_k}$ where $T,t(\tau)\subseteq\omega^*$ for $\tau\in T$. We associate with any  $\tau\in T$ the $k$-partition $\rho_T(\tau)$ by induction on the rank $rk(\tau)$ of $\tau$ in $(T;\sqsupseteq)$ as follows: if $rk(\tau)=0$ then $\rho_T(\tau)=\nu(V)$ where $V=t(\tau)\in\mathcal{T}_k$; otherwise, $\rho_T(\tau)=\nu(V)\cdot(\bigoplus\{\rho_T(\tau n)\mid n<\omega,\tau n\in T\})$. 
Finally, define a function $\rho:\mathcal{T}_{\mathcal{T}_k}\to k^{X^\omega}$ by $\rho(T)=\rho_T(\varepsilon)$. 

The main result of this subsection is the following.

\begin{Proposition}\label{wadred1}
The function $\rho$ induces an isomorphic embedding of the quotient-poset of $({\mathcal F}_{{\mathcal T}_k};\leq_h)$ into that of $(\mathcal{A}_k;\leq_{CA})$.
 \end{Proposition}

{\em Proof.} First we show that for all $T,V\in\mathcal{T}_{\mathcal{T}_k}$ we have: $\rho(T)\in I(k^{X^*};\leq_{CA},\oplus)$, and $T\leq_hV$ iff $\rho(T)\leq_{CA}\rho(V)$. This  easily follows from Proposition 16 in \cite{s17}. Note that the definitions of $\mu,\nu,\rho$ in \cite{s17} are for $k$-partitions of $\omega^\omega$ ($\omega$-words over $\omega$) rather than for $X^\omega$. The main difference between $\omega^\omega$ and $X^\omega$ is that the latter space in compact while the former one is not. This difference is essential only for the $\omega$-ary version of $\oplus$; since here we deal only with binary version of $\oplus$, the definitions and arguments of \cite{s17} work also here. Another  difference is that here we use some additional coding to preserve aperiodicity; this is also not essential up to $\equiv_{CA}$.

By Proposition \ref{slattice}, the embedding $\rho$ of $({\mathcal T}_{{\mathcal T}_k};\leq_h)$ into  $(\mathcal{A}_k;\leq_{CA})$ uniquely extends to a semilattice embedding (also denoted by $\rho$) of $({\mathcal F}_{{\mathcal T}_k};\leq_h)$ into  $(\mathcal{A}_k;\leq_{CA})$.
 \qed

The results above  imply the following.

\begin{Proposition}\label{wadred3}
For any $F\in{\mathcal F}_{{\mathcal T}_k}$ we have: $\rho(F)\in\mathcal{A}(F)$ and $\rho(F)\equiv_{CA}\bfSig(F)$. The FHs over the 2-bases $\mathcal{R}\bfSig$ and $\mathcal{A}\bfSig$ do not collapse.
 \end{Proposition}

{\em Proof.} The assertion $\rho(F)\in\mathcal{A}\bfSig(F)$ is checked by the usual induction on $F$, assuming it to be minimal. The assertion $\rho(F)\equiv_{CA}\bfSig(F)$ follows from the results in \cite{s17}. These assertions and the non-collapse of the FH over $\bfSig$ (that follows from the much more general results in \cite{km19}) implies the non-collapse of the FHs over $\mathcal{R}\bfSig$ and $\mathcal{A}\bfSig$.
 \qed

\subsection{Relating the FHs over $\mathcal{C}_\mathcal{M}$ and $\bfSig$}\label{relating}

To prove other properties of the embedding $\rho$, we have to establish close relations between FHs over $\bfSig$ and over a 2-base $\mathcal{C}_\mathcal{M}$ constructed from  a given automaton $\mathcal{M}$. For this we use preorders $\leq_0$ and $\leq_1$ on $C_\mathcal{M}$ defined in \cite{wag79} as follows: $c\leq_0d$, if some (equivalently, every) state in $d$ is reachable from some (equivalently, every) state in $c$; $c\leq_1d$, if $c\supseteq d$. The following is a reformulations of the corresponding facts observed in \cite{wag79}.

\begin{Lemma}\label{cicles}
 The structure $(C_\mathcal{M};\leq_0,\leq_1)$ is a compatible 2-preorder.
 \end{Lemma}

As in the proof of Proposition \ref{alttree}, instead of $(C_\mathcal{M};\leq_0,\leq_1,A)$, where $A:C_\mathcal{M}\to\bar{k}$, it is sometimes useful to consider its ``copy'' $(C_\mathcal{M}/_{\equiv_0};\leq_0,d)$ in $\mathcal{P}_{\mathcal{P}_Q}$) where $d:C_\mathcal{M}/_{\equiv_0}\to\mathcal{P}_Q$ is given by $d([c]_0)=([c]_0\leq_1,A|_{[c]_0})$. Note that the equivalence classes in $C_\mathcal{M}/_{\equiv_0}$ bijectively correspond to the reachable strongly connected components (SCCs) of (the graph of) $\mathcal{M}$ while the induced partial order $\leq_0$ indicates which SCCs ore reachable from other (smaller) ones.

The next lemma is also from \cite{wag79} (see also Section 7 in \cite{s98}). We  reproduce the  proof  of one direction for a further generalisation. 

\begin{Lemma}\label{cicles1}
 If $A\subseteq C_\mathcal{M}$ is $\leq_0$-up (resp. $\leq_1$-up) in $C_\mathcal{M}$  then $f^{-1}_\mathcal{M}(A)\in\bfSig^0_1$ (resp. $f^{-1}_\mathcal{M}(A)\in\bfSig^0_2$), otherwise $\bfPi^0_1\leq_{CA}f^{-1}_\mathcal{M}(A)$ (resp. $\bfPi^0_2\leq_{CA}f^{-1}_\mathcal{M}(A)$).
 \end{Lemma}

{\em Proof of one direction.} Let $A$ be not $\leq_0$-up, i.e. $c\in A\not\ni d$ for some $c\leq_0d$; we show that the standard  $\bfPi^0_1$-complete language $L$ of words that do not contain letter 1, is Wadge reducible to $f^{-1}_\mathcal{M}(A)$. By the definition of $f_\mathcal{M}$ and $\leq_0,$ there are $z,u,v,w\in X^*$ and $q\in c,r\in d$ such that $in\cdot z=q$, $q\cdot u=q$, $c$ is the set of states in the run $q\cdot u$, $q\cdot w=r$, $r\cdot v=r$, and $d$ is the set of states in the run $r\cdot v$. We define the (synchronous) continuous function $g$ on $X^\omega$ as follows. We scan subsequent letters of an input $x\in X^\omega$ waiting for the first occurrence of 1; while we do not see it, we construct $g(x)$ as $zuu\ldots$ (thus, if $x$ does not contain 1 at all then $g(x)=zu^\omega$, hence $f_\mathcal{M}(g(x))=c$); once we see the smallest $i$ with $x(i)=1$, we further construct $g(x)$ as $zu^{i+1}wv^\omega$,  hence $f_\mathcal{M}(g(x))=d$). Then $g$ $CS$-reduces $L$ to $f^{-1}_\mathcal{M}(A)$. 

Let now $A$ be not $\leq_1$-up, i.e. $c\in A\not\ni d$ for some $c\leq_1d$; we show that the standard  $\bfPi^0_2$-complete language $L$ of words that  contain infinitely many entries of 1, is Wadge reducible to $f^{-1}_\mathcal{M}(A)$. By the definition of $f_\mathcal{M}$ and $\leq_1,$ there are $z,u,v\in X^*$ and $q\in d\subseteq c$ such that $in\cdot z=q$, $q\cdot u=q$, $d$ is the set of states in the run $q\cdot u$, $q\cdot v=q$,  and $c$ is the set of states in the run $q\cdot v$. We define the (synchronous) continuous function $g$ on $X^\omega$ as follows. We again scan  $x$ looking for entries of 1; while there are no occurrences, we construct $g(x)$ as $zuu\ldots$ (thus, if $x$ does not contain 1 at all, $g(x)=zu^\omega$, hence $f_\mathcal{M}(g(x))=d$); if we see the smallest $i_0$ with $x(i_0)=1$, we continue to construct $g(x)$ as $zu^{i_0+1}vuu\ldots$ (thus, if there is precisely one entry then $g(x)=zu^{i_0+1}vu^\omega$, hence $f_\mathcal{M}(g(x))=d$); if we meet the second entry $x(i_1)=1,i_1>1$, we continue to construct $g(x)$ as $zu^{i_0+1}vu^{i_1+1}vuu\ldots$ (thus, if there are precisely 2 entries then $g(x)=zu^{i_0+1}vu^{i_1+1}vu^\omega$, hence $f_\mathcal{M}(g(x))=d$), and so on (thus, if there are infinitely many entries of 1 at positions $i_0<i_1<\cdots$ then our construction yields $g(x)=zu^{i_0+1}vu^{i_1+1}vu^{i_2+1}\cdots$, hence $f_\mathcal{M}(g(x))=c$). Thus, $g$ $CS$-reduces $L$ to $f^{-1}_\mathcal{M}(A)$.
 \qed

The next lemma generalizes the argument in the previous one. Let $\mathcal{C}_\mathcal{M}$ be the 2-base associated with $(C_\mathcal{M};\leq_0,\leq_1)$ as in Subsection \ref{sample}.

\begin{Lemma}\label{cicles3}
 For  every $F\in\mathcal{F}_{\mathcal{T}_k}$, if  $F\leq_h(C_\mathcal{M};\leq_0,\leq_1,A)$ then $\bfSig(F)\leq_{CA}A\circ f_\mathcal{M}$.
 \end{Lemma}

{\em Proof.} Let $\varphi$ be a witness for $F\leq_h(C_\mathcal{M};\leq_0,\leq_1,A)$ as in the proof of Theorem \ref{alttree}, i.e. a monotone function $\varphi:(F\sqsubseteq)\to(C_\mathcal{M};\leq_0)$ such that $t(\tau)\leq_hd([\varphi(\tau)])$ for each $\tau\in T$. For every $\tau\in F$, let $\varphi_\tau:(t(\tau);\sqsubseteq)\to([\varphi(\tau)]_0;\leq_1))$ be a monotone function such that $v(\sigma)=A(\varphi_\tau(\sigma))$ for each $\sigma\in t(\tau)=(V,v)$.

Let $q:C_\mathcal{M}\to Q$ be a function such that $q(c)\in c$ for every $c\in C_\mathcal{M}$. For every $c\in C_\mathcal{M}$, choose a word $v_{c}\in X^*$ such that $q(c)\cdot v_{c}=q(c)$ and the set of states in the run $q(c)\cdot v_{c}$ coincides with $c$. For every $c,d\in C_\mathcal{M}$ with $c\leq_0d$, choose a word $w_{c,d}\in X^*$ such that $q(c)\cdot w_{c,d}=q(d)$.

We have to Wadge reduce any $k$-partition $B\in\bfSig(F)$  to $A\circ f_\mathcal{M}$, assuming w.l.o.g. that the forest $F$ is minimal in the sense of Subsection \ref{lab}. Since the 2-base $\bfSig$ is reducible, by Proposition \ref{fhkp}(2) there is a reducible $F$-family $(\{U_\tau\},\{U_{\tau\sigma}\})$ that determines $B$. We now (somewhat informally) describe a $CA$-function $g$ (that is in fact a $CS$-function if $F$ is a tree) that reduces $B$ to $A\circ f_\mathcal{M}$. We argue by induction on the rank of $F$ in  $\mathcal{F}_{\mathcal{T}_k}$ and by induction on $|F|$. First we consider the case when $F$ is a tree (hence $U_\varepsilon=X^\omega$), and at the end of the proof we discuss corrections that should be made if $F$ is a proper forest. 

Let first  $F\in\mathcal{T}_{k}$, then the $F$-family that determines $B$ simplifies to $\{U_\tau\}$, where the components $\tilde{U}_\tau$, $\tau\in F$, partition $X^\omega$ to differences of $\bfSig^0_1$-sets, and  $f(\tau)=A(\varphi(\tau))$ whenever $x\in\tilde{U}_\tau$. Using the standard mind-change procedure for open sets $U_\tau$, it is straightforward to generalise the $\bfPi^0_1$-strategy in the proof of Lemma \ref{cicles1} (using the words $w_{c,d}$ for $c,d\in\varphi(F)$) and construct a $CS$-function $g$ such that $f_\mathcal{M}(x)=\varphi(\tau)$ whenever $x\in\tilde{U}_\tau$. Then $B\leq_{CS}A\circ f_\mathcal{M}$ via $g$, as desired.

Let now $F=s(T)$ for a minimal tree $T\in\mathcal{T}_{k}\setminus\bar{k}$, then the $F$-family that determines $B$ simplifies to $\{U_{\varepsilon\sigma}\}_{\sigma\in T}$, where the components $\tilde{U}_{\varepsilon\sigma}$, $\sigma\in f(\varepsilon)=(T,t)$, partition $X^\omega$ to differences of $\bfSig^0_2$-sets, and  $t(\sigma)=A(\varphi(\tau))$ whenever $x\in\tilde{U}_{\varepsilon\sigma}$. Using the standard mind-change procedure for $\bfSig^0_2$-sets $U_{\varepsilon\sigma}$, it is straightforward to generalise the $\bfPi^0_2$-strategy in the proof of Lemma \ref{cicles1} (using the words $v_{d}$ for $d\in[c]_0$, $c=\varphi(\varepsilon)$) and construct a $CS$-function $g$ such that $f_\mathcal{M}(x)=\varphi_\varepsilon(\sigma)$ whenever $x\in\tilde{U}_{\varepsilon\sigma}$. Then $B\leq_{CS}A\circ f_\mathcal{M}$ via $g$, as desired.

In the ``general'' case $F\in(\mathcal{T}_{k}(2)\setminus\mathcal{T}_{k})$, $F$ is non-singleton, we combine the constructions of the two previous paragraphs in the obvious way: while $x$ sits in $\tilde{U}_\varepsilon$, we follow the strategy of the previous paragraph working with the SCC $[\varphi(\varepsilon)]_0$ which gives us the reduction  $B\leq_{CS}A\circ f_\mathcal{M}$ restricted to $\tilde{U}_\varepsilon$; if $x$ moves to some $U_i$, $i\in\omega\cap F$, we discover this at a finite step (since all $U_\tau$ are open) and move (using the word $w_{\varepsilon i}$) to the SCC $[\varphi(i)]_0$; using the strategy of the previous paragraph for $\tilde{U}_i$ in place of $\tilde{U}_\varepsilon$ and $\varphi_i$ in place of $\varphi_\varepsilon$, gives us the reduction  $B\leq_{CS}A\circ f_\mathcal{M}$ restricted to $\tilde{U}_i$; and so on.

It remains to consider the case when $F=F_0\sqcup\cdots\sqcup F_n$, $n>0$, is a proper forest canonically decomposed to trees. Then $\{U_i\}_{i\in\omega\cap F}$ is a clopen partition of $X^\omega$, so for any $x$ there are $i,j$ such that $x\upharpoonright_j\cdot X^\omega\subseteq U_i$. Thus, given $x$, we can first find such $i,j$, then find $w_x\in X^*$ with $in\cdot w_x\in\varphi(i)$, and then follow the strategy of the previous paragraph with $U_i$ in place of $X^\omega$ and with $\varphi(i)$ in place of $\varphi(\varepsilon)$. In this way we obtain a $CA$-witness $g$, $g(x)\sqsupseteq w_x$, for $B\leq_{CA}A\circ f_\mathcal{M}$.
 \qed

The last lemma of this subsection extends Lemma \ref{cicles1} to all levels of the FH of $k$-partitions over  $\mathcal{C}_\mathcal{M}$.

\begin{Lemma}\label{cicles2}
 For all $A:C_\mathcal{M}\to\bar{k}$ and $F\in\mathcal{F}_{\mathcal{T}_k}$ we have: If $A\in\mathcal{C}_\mathcal{M}(F)$  then $A\circ f_\mathcal{M}\in\bfSig(F)$, otherwise $\bfSig(T)\leq_{CA}A\circ f_\mathcal{M}$ for some $T\in M(F)$ and hence $A\circ f_\mathcal{M}\not\in\bfSig(F)$.
 \end{Lemma}

{\em Proof.} By Lemma \ref{cicles1}, $f^{-1}_\mathcal{M}:\mathcal{C}_\mathcal{M}\to\bfSig$ is a morphism of the 2-bases, hence $A\in\mathcal{C}_\mathcal{M}(F)$  then $A\circ f_\mathcal{M}\in\bfSig(F)$ by Proposition \ref{fhkp}(4) which yields the first assertion. Let now $A\not\in\mathcal{C}_\mathcal{M}(F)$. By Theorem \ref{alttree}, $T\leq_h(C_\mathcal{M};\leq_0,\leq_1,A)$ for some $T\in M(F)$. By Lemma \ref{cicles3}, $\bfSig(T)\leq_{CA}A\circ f_\mathcal{M}$. Since $T\not\leq_hF$, $\rho(T)\not\leq_{CA}\rho(F)$ by Proposition \ref{wadred1}, hence $A\circ f_\mathcal{M}\not\in\bfSig(F)$ by Proposition \ref{wadred3}. 
 \qed

\subsection{Main results}\label{wadged}

Here we complete the topological classification of (aperiodic) regular languages by showing that the embedding $\rho$ is surjective, and then deduce a series of corollaries. 

The next result together with Proposition \ref{wadred1} imply Theorem \ref{wadred}.

\begin{Proposition}\label{wadred2}
For every Muller $k$-acceptor $(\mathcal{M},A)$ there is $F\in{\mathcal F}_{{\mathcal T}_k}$ such that  $A\circ f_\mathcal{M}\equiv_{CA}\rho(F)$.
 \end{Proposition}

{\em Proof.} Let $(F,f)$ be the forest unfolding of $(C_\mathcal{M};\leq_0,\leq_1,A)$ (more precisely, again of its copy in ${\mathcal P}_{{\mathcal P}_k}$). Since every equivalence class in $C_\mathcal{M}/_{\equiv_0}$ has the least element w.r.t. $\leq_1$ by Lemma \ref{cicles}, we have $F\in{\mathcal F}_{{\mathcal T}_k}$ by the remarks in Subsection \ref{lab}. It remains to show that $\bfSig(F)\leq_{CA}A\circ f_\mathcal{M}\in\bfSig(F)$. The reduction $\bfSig(F)\leq_{CA}A\circ f_\mathcal{M}$ holds by Lemma \ref{cicles3} because $f$ is a witness for $T\leq_h(C_\mathcal{M};\leq_0,\leq_1,A)$. To prove that $A\circ f_\mathcal{M}\in\bfSig(F)$, it suffices by Lemma \ref{cicles2} to show that $A\in\mathcal{C}_\mathcal{M}(F)$.  By Theorem \ref{alttree}, it suffices to show that if $T\leq_h(C_\mathcal{M};\leq_0,\leq_1,A)$ and $T\in\mathcal{T}_{\mathcal{T}_k}$ then $T\leq_hF$. This follows from the fact mentioned in Subsection \ref{lab} that $F$ is a largest forest below $(C_\mathcal{M};\leq_0,\leq_1,A)$ w.r.t. $\leq_h$.
 \qed

The next assertion is an immediate corollary of Theorems \ref{wadred} and \ref{fotheory}.

\begin{Corollary}\label{fotheo}
If $k=2$ then the first-order theories of the quotient-posets of $(\mathcal{R}_k;\leq_{CA})$ and $(\mathcal{A}_k;\leq_{CA})$ are decidable, otherwise they are undecidable and, moreover, computably isomorphic to the first-order arithmetic.  
 \end{Corollary}

Now we slightly strengthen Theorem \ref{wadred} by expanding the signature. Let $\mathbb{R}_k$ be the quotient-structure   of $(\mathcal{R}_k;\leq_{CA},I,\oplus,\cdot,q_0,\ldots,q_{k-1})$ under  $\equiv_{CA}$ where $I$ is the unary relation true precisely on the join-irreducible elements, and $\cdot,q_0,\ldots,q_{k-1}$ are the restrictions of the operations from Subsection \ref{reduc0} to $\mathcal{R}_k$. Let the structure $\mathbb{A}_k$ be defined similarly but on the universe $\mathcal{A}_k$. Let $\mathbb{F}_{\mathcal{T}_k}$ be the quotient-structure   of $(\mathcal{F}_{\mathcal{T}_k};\leq_{h},I,\sqcup,\cdot,q_0,\ldots,q_{k-1})$ under  $\equiv_{h}$ where $I$ is defined similarly, and $q_i(F)=s(i\cdot r(F))$ for every $F\in\mathcal{F}_{\mathcal{T}_k}$.

\begin{Proposition}\label{wadred4}
The structures $\mathbb{R}_k$, $\mathbb{A}_k$, and $\mathbb{F}_{\mathcal{T}_k}$ are isomorphic.
 \end{Proposition}

{\em Proof Sketch.} The  isomorphisms are of course induced by the embedding $\rho$. By Theorem \ref{wadred}, this  induced function is an isomorphism of the quotient-orders. Since the relation $I$ and the operation of supremum are defined in terms of the ordering, they are also preserved by the function. It remains to show that the function also respects the operation $\cdot,q_0,\ldots,q_{k-1}$. By induction on the cardinality of the involved forests it is not hard to check that $\rho(F\cdot G)\equiv_{CA}\rho(F)\cdot\rho(G)$ and $\rho(q_i(F))\equiv_{CA}q_i(\rho(F))$.
 \qed

We have completed the classifications of (aperiodic) regular sets up to the Wadge reducibility $\leq_{CA}$ but this reducibility is not well suited for automata theory because it is highly non-effective. That was the main reason to consider the $DA$-reducibility on $\mathcal{R}_k$ \cite{wag79} and the $AA$-reducibility on $\mathcal{A}_k$ \cite{s08}. We show that these reducibilities behave on the corresponding ``automatic'' FHs in the right way in the sense of hierarchy theory \cite{s08a}.

\begin{Proposition}\label{wadred5}
$DA$-Reducibility (resp. $AA$-reducibility) fits the FH over $\mathcal{R}\bfSig$ (resp. $\mathcal{A}\bfSig$).
 \end{Proposition}

{\em Proof.} ``Fits'' in the formulation means that every level of the FH over $\mathcal{R}\bfSig$ (resp. $\mathcal{A}\bfSig$) is closed under $\leq_{DA}$ (resp. $\leq_{AA}$). Every level of the Borel hierarchy is closed under both $\leq_{DA},\leq_{AA}$, the class $\mathcal{R}$ is closed under $\leq_{DA}$ \cite{wag79}, and $\mathcal{A}$ is closed under $\leq_{AA}$. Thus, both $\mathcal{R}\cap\bfSig^0_1$ and $\mathcal{R}\cap\bfSig^0_2$ are closed under $\leq_{DA}$, while both $\mathcal{A}\cap\bfSig^0_1$ and $\mathcal{A}\cap\bfSig^0_2$ are closed under $\leq_{AA}$. Thus, for every $DA$-function (resp. $AA$-function) $f$ the preimage map $f^{-1}$ is a morphism of the 2-base $\mathcal{R}\bfSig$ (resp. $\mathcal{A}\bfSig$) into itself. By Proposition \ref{fhkp}(4), every level $\mathcal{R}\bfSig(F)$ (resp.  $\mathcal{R}\bfSig(F)$) is closed under every such $f^{-1}$, i.e. it is  closed under $\leq_{DA}$ (resp. $\leq_{AA}$). 
 \qed

Next we prove the analogue of Proposition \ref{redcoin1} for $CA$-reducibility. Together with Proposition \ref{redcoin1}, the next theorem extends to $k$-partitions the corresponding important facts about sets proved in \cite{wag79,s08}.

\begin{Theorem}\label{redcoin}
 The relation $\leq_{CA}$  coincides with $\leq_{DA}$ on $\mathcal{R}_k$, and with $\leq_{AA}$ on $\mathcal{A}_k$.
\end{Theorem}

{\em Proof.} Both assertions are proved similarly, so we prove only the second one. We have to show that, for every $A,B\in\mathcal{A}_k$, $A\leq_{CA}B$ implies $A\leq_{AA}B$. By Theorem \ref{wadged}, $A\equiv_{CA}\rho(F)$ and $B\equiv_{CA}\rho(G)$ for some (minimal) forests $F,G$. If $G$ is a tree then $\rho(F)\leq_{CS}\rho(G)$ by the proof of Lemma \ref{cicles3}, hence $\rho(F)\leq_{AS}\rho(G)$ by Proposition  \ref{redcoin1}, hence $\rho(F)\leq_{AA}\rho(G)$.

Let now $G=G_0\sqcup\cdots\sqcup G_n$, $n>0$, be a proper forest canonically decomposed to trees. If $F$ is a tree then $F\leq_hG$ by Theorem \ref{wadred}, hence $F\leq_hG_j$ for some $j\leq n$. Then $\rho(F)\leq_{CA}\rho(G_j)$ by Theorem \ref{wadred}, hence  $\rho(F)\leq_{AA}\rho(G_j)\leq_{AA}\rho(G)$ by the previous paragraph. Finally, let also $F=F_0\sqcup\cdots\sqcup F_m$, $m>0$, be a proper forest canonically decomposed to trees. Then $F_i\leq_hG$ for all $i\leq m$, hence $\rho(F_i)\leq_{CA}\rho(G)$ for all $i\leq m$. By the previous case, $\rho(F_i)\leq_{AA}\rho(G)$ for all $i\leq m$, hence $\rho(F)\equiv_{AA}\rho(F_0)\oplus\cdots\oplus\rho(F_m)\leq_{AA}\rho(G)$.
 \qed
 
Next we prove an important preciseness property of the Boolean algebras $\mathcal{R}$ and $\mathcal{A}$ (see the end of Subsection \ref{wadged}). In general, precise Boolean algebras are rare but for these ones we have the following. 

\begin{Theorem}\label{wadred6}
The Boolean algebras $\mathcal{R}$ and $\mathcal{A}$ are $\bfSig$-precise. 
 \end{Theorem}

{\em Proof.} The proofs for both algebras are similar, so we consider only $\mathcal{A}$. By the definition at the end of Subsection \ref{wadged}, we have to show that $\mathcal{A}_k\cap\bfSig(F)\subseteq\mathcal{A}\bfSig(F)$ for every $F\in{\mathcal F}_{{\mathcal T}_k}$. Let $B\in\mathcal{A}_k\cap\bfSig(F)$, then $B\equiv_{CA}\rho(G)$ for some $G\leq_hF$ by Theorem \ref{wadred}. Since both $B$ and $\rho(G)$ are in $\mathcal{A}_k$, we get $B\equiv_{AA}\rho(G)$ by Theorem \ref{redcoin}. It suffices to show that $\rho(G)$ is in $\mathcal{A}\bfSig(G)$ (because then also $B$ is in $\mathcal{A}\bfSig(G)\subseteq\mathcal{A}\bfSig(F)$ by Propositions \ref{wadred5} and \ref{fhkp}(1)). The assertion $\rho(G)\in\mathcal{A}\bfSig(G)$ is checked by a straightforward induction on the size of the forest $G$ (that may w.l.o.g. be assumed minimal), using the definition of $\rho$ in Subsection \ref{reduc0} and the fact that the 2-base $\mathcal{A}\bfSig(G)$ is reducible \cite{s08}.
 \qed

The next proposition gives a very clear description of the relationship between the FHs over $\mathcal{R}\bfSig$ and over the 2-bases $\mathcal{C}_\mathcal{M}$.

\begin{Proposition}\label{wadred7}
For every  $F\in{\mathcal F}_{{\mathcal T}_k}$, we have: $\mathcal{R}\bfSig(F)=\bigcup\{A\circ f_\mathcal{M}\mid A\in\mathcal{C}_\mathcal{M}(F)\}$.
 \end{Proposition}

{\em Proof.} The inclusion $\supseteq$ follows from Lemma \ref{cicles2} and Theorem \ref{wadred6}. Conversely, let $B\in\mathcal{R}\bfSig(F)$, then $B\in\mathcal{R}_k$ and $B\in\bfSig(F)$ by Proposition \ref{fhkp}(8), hence $B=A\circ f_\mathcal{M}$ for some Muller $k$-acceptor $(\mathcal{M},A)$. It suffices to show that $A\in\mathcal{C}_\mathcal{M}(F)$. Suppose not, then by Theorem \ref{alttree}, there is $T\in\mathcal{T}_{\mathcal{T}_k}$ such that $T\leq_h(C_\mathcal{M};\leq_0,\leq_1,A)$ and $T\not\leq_hF$. By Lemma \ref{cicles3}, $\bfSig(T)\leq_{CA}A\circ f_\mathcal{M}$. Thus, $\bfSig(T)\leq_{CA}B\in\bfSig(F)$, so $\bfSig(T)\subseteq\bfSig(F)$, contradicting Theorem \ref{wadred}.
 \qed

%\subsection{The structures of Lipschitz degrees}\label{lip}

 We conclude this subsection by remarks on the structure of $CS$-degrees of (aperiodic) regular sets (as already mentioned, in descriptive set theory they are known as Lipschitz degrees).
The next assertion provides some information relating the Wadge and Lipschitz reducibilities. We use the usual notation about forests and trees mentioned above.

\begin{Proposition}\label{relfh}
\begin{enumerate}\itemsep-1mm
\item  If at least one of the minimal forests $F,G$ is a tree then $\rho(F)\leq_{CA}\rho(G)$ iff $\rho(F)\leq_{CS}\rho(G)$.
 \item If $F$ is a tree then $[\rho(F)]_{CA}=[\rho(F)]_{CS}$, i.e. the Wadge degree of $\rho(F)$ consists of a single Lipschitz degree, otherwise it splits to infinitely many Lipschitz degrees.
  \end{enumerate}
 \end{Proposition}

{\em Proof.} (1) Follows from the proof of Theorem \ref{redcoin}.

(2) If $F$ is a tree, the assertion follows from (1), so let $F=F_0\sqcup\cdots\sqcup F_m$, $m>0$, be a proper forest canonically decomposed to trees. Then $[\rho(F)]_{CA}$ contains an $\omega$-chain of Lipschitz degrees that is constructed  in the same way as in the particular case of 2-partitions (see \cite{wag79}, Theorem 9.1 in \cite{s98} and Theorem 10 in \cite{s08}).
 \qed

In fact, it is not  hard (but a bit cumbersome) to characterise the quotient-posets of $(\mathcal{R}_k;\leq_{CS})$ and of $(\mathcal{A}_k;\leq_{CS})$ up to isomorphism (e.g., if  $m=1$ then $[\rho(F)]_{CA}$ splits to an $\omega$-chain of Lipschitz degrees, but for $m>1$ the structure of Lipschitz degrees inside $[\rho(F)]_{CA}$ becomes more complicated but remains understandable). These characterisations imply that
the quotient-posets of $(\mathcal{R}_k;\leq_{CS})$ and of $(\mathcal{A}_k;\leq_{CS})$ are isomorphic. This extends to $k$-partitions the corresponding fact for sets implicitly contained in \cite{wag79,s09}.

\section{Computability and complexity issues}\label{compl}

There are many natural algorithmic problems related to topological properties  of regular  sets considered e.g. in \cite{wag79,kpb,wy}. Here we briefly discuss extensions of these problems to $k$-partitions, and  some new algorithmic problems.

First we discuss algorithmic problems that apparently were  not considered explicitly  in the literature on automata theory but are very popular in computability theory where people are interested in characterizing the complexity of presentation of natural countably infinite algebraic structures of finite signatures. Such a structure is {\em computably presentable}  if it is isomorphic to a structure whose universe is $\omega$  and all signature functions and relations are computable.  

In preceding subsections we  considered several natural structures including the structures $\mathbb{R}_k$ and $\mathbb{A}_k$ from Proposition \ref{wadred4}. Note that from their definition it is hard to see that they are computably presentable.

\begin{Proposition}\label{compresent}
The structures $\mathbb{R}_k$ and $\mathbb{A}_k$ are computably presentable.
 \end{Proposition}

{\em Proof.} By Proposition \ref{wadred4}, bot structures are isomorphic to $\mathbb{F}_{\mathcal{T}_k}$, hence it suffices to show that the latter structure is computably presentable. Considering only ``concrete'' trees and forests, and remembering the definitions in Subsection \ref{lab}, we see that there is a natural effective surjection (naming) $\nu$ from $\omega$ onto $\mathcal{F}_{\mathcal{T}_k}$ such that $\leq_h$, $r$, and $\sqcup$ are represented by computable relations and functions on the names (so, e.g., the relation $\nu(m)\leq_h\nu(n)$ is computable). Moreover, there is a computable function $f$ that finds, given any $n$, a minimal forest $\nu(f(n))$ $h$-equivalent to $\nu(n)$. This implies that the relation $I(\nu(n))$ is computable (because it is equivalent to $\nu(f(n))$ being a tree). In the same manner we check that the functions $q_0,\ldots,q_{k-1},\cdot$ are represented by computable functions on the names. From standard facts of computability theory it now follows that $\mathbb{F}_{\mathcal{T}_k}$ is computably presentable.
 \qed

From the remarks at the end of the previous section it follows that the structures of Lipschitz degrees of (aperiodic) regular $k$-partitions (i.e., the quotient-posets of $(\mathcal{R}_k;\leq_{DS})$ and $(\mathcal{A}_k;\leq_{AS})$) are also computable presentable.

The problem of finding feasible presentations of a given structure is more subtle.
A structure is {\em $p$-presentable} if there is a surjection from a polynomial-time computable subset of $X^*$ onto the universe of the structure modulo which all signature  functions and relations, and also the equality relation, are polynomial-time computable. We abbreviate ``polynomial-time computable'' to ``$p$-computable''.

Jointly with P.E. Alaev we have recently shown that the structures in Proposition \ref{compresent} are in fact $p$-presentable. This result (among others) should be published separately. The proof extends the coding and proofs  in \cite{hs14} where some particular cases are considered. 

The computational complexity of functions and relations about regular languages are usually studied when the languages are given by their standard ``names'' like automata or regular expressions. In our context it is natural to think that $k$-partitions are given by Muller's $k$-acceptors recognising them. In particular, for the relation $I$ one could wish to estimate the complexity of the problem: given an aperiodic Muller's $k$-acceptor $(\mathcal{M};A)$, is the $k$-partition $A\circ f_\mathcal{M}$ join-irreducible in $(\mathcal{A}_k;\leq_{CA},\oplus)$? For the function $\cdot$, one could be interested in estimating the complexity of the problem: given  aperiodic Muller's $k$-acceptors for $A,B\in\mathcal{A}_k$, find an aperiodic Muller's $k$-acceptor for $A\cdot B$ (up to $\equiv_{CA}$). From the results above we easily obtain the following.

\begin{Corollary}\label{compu}
All the signature functions and relations on $\mathcal{A}_k$ in Proposition \ref{compresent} are computable w.r.t. the Muller $k$-acceptor presentation.
\end{Corollary}

{\em Proof.} Consider e.g. the relation $\leq_{CA}$. Given acceptors $(\mathcal{M};A)$ and $(\mathcal{M}_1;A_1)$ recognising resp. $A\circ f_\mathcal{M}$ and $A_1\circ f_{\mathcal{M}_1}$, compute, using the algorithms in the proofs of the corresponding facts above, (names of) $F,F_1\in\mathcal{F}_{\mathcal{T}_k}$ such that $A\equiv_{CA}\rho(F)$ and $A_1\equiv_{CA}\rho(F_1)$, and check $F\leq_hF_1$, using the computable presentation in the proof of Proposition \ref{compresent}.
 \qed

The method of Corollary \ref{compu} and the computability of many other relations and functions on the wqo $\mathcal{F}_{\mathcal{T}_k}$ imply the computability of many other topological problems about regular (aperiodic) $k$-partitions. The complexity of such problems is much more subtle and leads to  interesting open questions. Even the $p$-computability of $\leq_{CA}$ for $k>2$ is currently open because our approach involves  computing of the forest unfolding of a $\bar{k}$-labeled 2-preorder  (see Subsection \ref{lab}) in polynomial time which is easy for $k=2$ but far from obvious for $k>2$. For $k=2$, the $p$-computability of $\leq_{CA}$ and of some other related functions and relations is known from \cite{kpb,wy}.

\section{Other classes of $k$-partitions of $\omega$-words}\label{other}

In this section we briefly discuss some other classes of $k$-partitions of $\omega$-words for which the method described above could help.

 We look at some classes $\mathcal{D}$ of $\omega$-languages recognised by relatively simple computing devices (the class of context-free $\omega$-languages is considered as too wide, according to the non-decidability results for this class mentioned in the Introduction). We briefly discuss  some such classes $\mathcal{D}$ divided into 3 categories: those   in between $\mathcal{A}$ and $\mathcal{R}$, those  below $\mathcal{A}$, and those beyond $\mathcal{R}$.

For classes $\mathcal{D}$ with $\mathcal{A}\subset\mathcal{D}\subset\mathcal{R}$, there is no problem with characterising the corresponding Wadge and Lipschitz  degrees since, by Theorem \ref{wadred} and the remarks at the end of Subsection \ref{wadged}, the quotient-posets of $(\mathcal{D}_k;\leq_{CA})$ and $(\mathcal{D}_k;\leq_{CS})$ are isomorphic respectively to $(\mathcal{R}_k;\leq_{CA})$ and $(\mathcal{R}_k;\leq_{CS})$. But analogues of other results above are non-trivial and interesting.

We did not find (excepting \cite{rt}) in the literature papers investigating classes  $\mathcal{D}$ with $\mathcal{A}\subset\mathcal{D}\subset\mathcal{R}$, but there are investigations of classes of regular languages of finite words containing the regular aperiodic languages (see e.g. \cite{str,s09}). One could consider classes of $\omega$-languages obtained in a way similar to the definition of the aperiodic $\omega$-languages from the aperiodic languages of finite words. We apply this idea to the classes of quasi-aperiodic and of $d$-aperiodic languages (for every  integer $d>1$) \cite{str,s09}.

Associate with any alphabet $X$ the signature
$\sigma=\{\leq,Q_a\mid a\in X\}$ where $\leq$
is (a name of) the binary relation  interpreted as the usual order on positions of a word, and $Q_a$ is the unary
relation true at the positions of the letter $a$. By a theorem of McNaughton and Papert,  the class of
${\rm FO}_\sigma$-axiomatizable languages (i.e.,  languages satisfying a fixed first-order sentence
of $\sigma$), coincides with the class of regular aperiodic
languages; this class of languages also coincides with the class of languages recognised by Muller's acceptors without counting pattern (an automaton $\mathcal{M}$ {\em has a counting pattern} if there are $n>1$, a reachable state $q$, and a word $v\in X^*$ such that $q\cdot v^n=q$ and $q\cdot v^m\not=q$ for $m<n$). Similar facts  hold for $\omega$-languages.  

Let $\tau_d=\sigma\cup\{P_d\}$, where $P_d$ is the unary
relation  true on the positions of  a word which are divisible by $d$, and let $\tau=\bigcup_d\tau_d$. Then the ${\rm FO}_{\tau_d}$-axiomatizable (resp. ${\rm FO}_{\tau}$-axiomatizable) languages of non-empty finite words coincide with the so called $d$-aperiodic (resp. quasi-aperiodic) languages  \cite{str,s09}. Also, the ${\rm FO}_{\tau_d}$-axiomatizable (resp. ${\rm FO}_{\tau}$-axiomatizable) languages of finite non-empty words coincide with the languages recognised by $d$-aperiodic (resp. balanced-aperiodic) acceptors defined as follows. An automaton $\mathcal{M}$ {\em has a $d$-counting pattern} if there are $n>1$, a reachable state $q$, and a word $v\in X^*$ such that $q\cdot v^n=q$, $q\cdot v^m\not=q$ for $m<n$, and $d$ divides $|v|$. An automaton $\mathcal{M}$ {\em has a balanced counting pattern} if there are $n>1$, a reachable state $q$, and  words $u,v\in X^*$ such that $q\cdot v^n=q$, $q\cdot v^m\not=q$ and $(q\cdot v^m)\cdot u=q\cdot v^m$ for $m<n$, and $|u|=|v|$. An automaton is {\em $d$-aperiodic (resp. balanced-aperiodic)} if it has no $d$-counting (resp. balanced counting) patterns.  

We do not currently know whether the results of the previous paragraph hold for $\omega$-languages but it seems quite natural to take the introduced classes of automata  in place of aperiodic automata in the attempt to develop analogues of the above theory for $\mathcal{A}$ for the classes $\mathcal{A}^d,\mathcal{A}^\tau$ defined as follows. Let $\mathcal{A}^d$ (resp. $\mathcal{A}^\tau$) be the class of $\omega$-languages recognised by the $d$-aperiodic (resp. balanced aperiodic) Muller acceptors, then $\mathcal{A}\subset\mathcal{A}^d\subset\mathcal{A}^\tau\subset\mathcal{R}$. The {\em $dA$-functions (resp. $\tau A$-functions)} are those computed by the $d$-aperiodic (resp. balanced aperiodic) asynchronous transducers. The  $dS$-functions (resp. $\tau S$-functions) are defined similarly but with synchronous transducers. The corresponding reducibilities $\leq_{dA},\leq_{\tau A},\leq_{dS},\leq_{\tau S}$ are defined in the obvious way.  We guess that the  introduced notions have properties similar to those in Section \ref{prel} for aperiodic sets (in particular, we guess that the corresponding versions of the B\"uchi-Landweber theorem hold). If this is really the case, there should be no problem to extend the whole theory of this paper to the classes $\mathcal{A}^d_k$ and $\mathcal{A}^\tau_k$ of $k$-partitions.

Another interesting class  $\mathcal{A}^{MOD}$ of $\omega$-languages axiomatized by  $\sigma$-sentences with first-order and modulo quantifiers was considered in \cite{rt} (obviously, $\mathcal{A}^\tau\subset\mathcal{A}^{MOD}\subset\mathcal{R}$). It was shown that for this class (and the class of so called causal functions defined by $\{<,MOD\}$-formulas) the analogue of B\"chi-Landweber's theorem holds. It gives a hope that also analogues of the results of this paper hold for the corresponding class of $k$-partitions $\mathcal{A}^{MOD}_k$. For unification purposes, it seems natural to characterize the class of causal  $\{<,MOD\}$-functions in terms of suitable automata similar to those discussed above\footnote{I thank Wolfgang Thomas for a recent discussion of this problem.}.

In \cite{ch12} some natural subclasses of $\mathcal{A}$ (essentially, induced by some levels of the Brzozowski and Straubing-Th\'erien hierarchies (see e.g. \cite{pp,str,s09}) were investigated in a search for analogues of the B\"uchi-Landweber theorem. The results in \cite{ch12} suggest that for such classes the Wadge degrees may be characterised without big problems but e.g. natural reducibilities that fit the corresponding effective FHs (analogues of Proposition \ref{wadred5}) are hard to find (if they exist at all).

Among many natural superclasses of $\mathcal{R}$, the class of visibly push down (VPD) $\omega$-languages \cite{vonBraunmuehl_Verbeek,AlurMadhusudan} and its subclasses seem especially interesting. Such an investigation was recently initiated in \cite{okh} but the most inetresting questions remain open. For the deterministic VPD $\omega$-languages we expect that many results of this paper remain true but with the algorithmic problems some surprises are possible. For the non-deterministic VPD $\omega$-languages, there is still a hope to build a similar theory (due to a nice determinisation theorem established in \cite{LoedingMadhusudanSerre}), though the set of Wadge degrees occupied by such languages is certainly larger than the Wadge degrees of regular $\omega$-languages, as it follows from the results in
\cite{LoedingMadhusudanSerre}.

\section{Future work} \label{open}

We believe that the methods developed in this paper could help in realising the project sketched in the previous subsection. This project seems interesting not only for automata theory (because it sketches simple and potentially useful classifications of several natural classes of $k$-partitions) but also for descriptive set theory, as a first step in identification the constructive content of the Wadge theory.

\end{document}